\documentclass[12pt,a4paper]{article}
\usepackage{amssymb,amsmath,amsthm}
\usepackage[english]{babel}
\usepackage{t1enc}
\usepackage[latin2]{inputenc}
\usepackage{epsfig}
\usepackage{subfig}
\usepackage{epstopdf}
\usepackage{xcolor}
\usepackage{comment}
\usepackage{hyperref}
\usepackage{tabularx}
\usepackage[bottom=1.5in, top=1.25in]{geometry}

\newtheorem{thm}{Theorem}
\newtheorem{cor}[thm]{Corollary}
\newtheorem{defi}[thm]{Definition}
\newtheorem{lem}[thm]{Lemma}

\newtheorem{claim}[thm]{Claim}
\newtheorem{conj}[thm]{Conjecture}
\newtheorem{prob}[thm]{Problem}
\newtheorem{rem}[thm]{Remark}
\newtheorem{obs}[thm]{Observation}

\theoremstyle{remark}
\newtheorem{remark}[thm]{Remark}

\usepackage{indentfirst}
\def\01{0\text{-}1}

\makeatletter
\def\blfootnote{\gdef\@thefnmark{}\@footnotetext}
\makeatother

\title{Saturation of edge-ordered graphs}
\begin{document}
	\author{Vladimir Bo\v skovi\'c\thanks{Universit\'e Paris-Saclay, CNRS, CEA, Institut de Physique Th\'eorique, 91191 Gif-sur-Yvette, France\\Email address: vladimir.boskovic@ipht.fr} \and Bal\'azs Keszegh\thanks{HUN-REN Alfréd Rényi Institute of Mathematics and ELTE Eötvös Loránd University, Budapest, Hungary.\\Email address: keszegh@renyi.hu\\			}}
	\maketitle	
	
	\begin{abstract}
		For an edge-ordered graph $G$, an $n$-vertex edge-ordered graph $H$ is $G$-saturated if it is $G$-free and adding any new edge with an arbitrary label to $H$ creates a copy of $G$. The saturation function is the minimum number of edges in a $G$-saturated graph.
		
		For (unordered) graphs, $0$-$1$ matrices, and vertex-ordered graphs, the saturation function is always $O(n)$ and satisfies a dichotomy: it is either $O(1)$ or $\Theta(n)$. The saturation function of an edge-ordered graph follows a weaker dichotomy, being either $O(1)$ or $\Omega(n)$. However, by finding edge-ordered graphs whose saturation functions are $\Omega(n \sqrt{\log n})$, we show that $O(n)$ is not a universal upper bound. 
		
		We also study the semisaturation problem for edge-ordered graphs, a variant of the saturation problem in which $H$ is not required to be $G$-free. We prove a general upper bound $O(n \log n)$ and characterize edge-ordered graphs with bounded semisaturation functions. 
		
		We then present several families of edge-ordered graphs with bounded, linear, and superlinear (semi)saturation functions. We also introduce a natural variant of saturation in which the added edge is required to receive the smallest label. The behaviour of the two variants is similar in many respects, which motivated us to investigate the second variant extensively.
	\end{abstract}

\section{Introduction}

An \emph{edge-ordered graph} is a finite simple graph $G$ with a linear order on its edge set $E$. We usually assign the edge-order using an injective labeling $l: E \rightarrow \mathbb{N}$, but it can equivalently be any other linear order, e.g., $l: E \rightarrow \mathbb{R}$. We say that an edge-ordered graph $H$ \emph{contains} an edge-ordered graph $G$ if $G$ is isomorphic to a subgraph of $H$. Note that an isomorphism between two edge-ordered graphs has to respect the edge-order. If $H$ does not contain $G$, then we say that $H$ is $G$-free, or equivalently that $H$ \emph{avoids} $G$. 

A classical saturation problem of (not edge-ordered) graphs asks for the minimum number of edges in an $n$-vertex graph $H$ such that it avoids $G$ and adding any missing edge to $H$, a copy of $G$ is created. In general, we always assume that $n$ is a large enough integer and we denote by $sat(n, G)$, the saturation number, the number of edges in such a graph $H$. 

In the extremal problem of edge-ordered graphs~\cite{gerbner}, we are looking for the maximum size of the saturating host graph, in the sense that no edge can be added with any label.

\begin{defi} \label{ex_function}
	Let $G$ be an edge-ordered graph and $n \in \mathbb{N}$.		
	We define $ex_e(n,G)$ to be the maximum number of edges in a $G$-free edge-ordered graph $H$ on $n$ vertices.
\end{defi}

We extend the notion of saturation to edge-ordered graphs. However, when a missing edge is added, there are several possible ways to extend the linear order on the edges of $H$ by assigning a label to the new edge. Unlike vertex-ordered graphs that were studied in~\cite{boskesz023}, this is not uniquely defined. Therefore, we introduce three different variants.

\begin{defi} \label{definitionsAMS}
	Let $G$ be an edge-ordered graph and $n \in \mathbb{N}$.
	
	We define $sat_e(n,G)$ to be the minimum number of edges in a $G$-free edge-ordered graph $H$ on $n$ vertices such that adding any new edge with an arbitrary label (inserted anywhere in the linear order of edges) introduces a copy of $G$ (that is, every label is forbidden; $e$ stands for `every' or `edge-ordered').
	
	We define $sat_m(n,G)$ to be the minimum number of edges in a $G$-free edge-ordered graph $H$ on $n$ vertices in which adding any new edge with a minimal label (inserted at the beginning of the linear order of edges) introduces a copy of $G$ ($m$ stands for `minimal').
	
	We define $sat_s(n,G)$ to be the minimum number of edges in a $G$-free edge-ordered graph $H$ on $n$ vertices in which for any new edge there exists some label (inserted somewhere in the linear order of edges) for which adding this edge introduces a copy of $G$ ($s$ stands for `some').
	
	In each case, any (not necessarily minimal) graph $H$ satisfying the corresponding property is called a \emph{host graph}, and we say that $H$ saturates $G$. \footnote{It will be clear from the context which definition is being considered.}
\end{defi}

From the definitions, it follows every host graph for $sat_e$ is also a host graph for $sat_m$.  Similarly, every host graph for $sat_m$ is also a host graph for $sat_s$. Consequently, we obtain the following inequalities.

\begin{rem}\label{rem:sme}
	For every edge-ordered graph $G$, we have \[sat_s(n,G) \le sat_m(n,G) \le sat_e(n,G)\le ex_e(n,G).\]
\end{rem}

Therefore, the most natural saturation function is $sat_e(n,G)$, obtained from the extremal problem by searching for a minimal saturating host graph instead of a maximal one.\footnote{This is why $e$ refers to `edge-ordered' along `every'.} However, in many cases, instead of dealing with every possible label of the new edge, we already get the same behavior if we forbid the minimal and maximal labels. In fact, we already gain a lot of insight by looking at the minimal label only, which is slightly more convenient to handle. Therefore, our center of attention is $sat_m$ throughout the paper. The consequences for $sat_e$ are discussed only afterwards, while we consider $sat_s$ only briefly compared to the other two functions. 

We also consider a natural variant of the saturation problem in which $H$ is not assumed to be $G$-free. Therefore, we only ask for a graph with a minimum number of edges such that, by adding a new edge with every/minimal/some label, a new copy of $G$ is created. This version is called semisaturation, and the respective functions are denoted by $ssat_e/ssat_m/ssat_s$. By definition, the following properties hold and will be used throughout the paper.

\begin{rem}
	For every edge-ordered graph $G$, we have \[ssat_e(n,G) \le sat_e(n, G),~~~~
	ssat_m(n,G) \le sat_m(n,G),~~~~
	ssat_s(n,G) \le sat_s(n,G).\]
\end{rem}

\subsection{History}

Saturation problems were first studied by Erd\H os, Hajnal, and Moon~\cite{erdoshajnalmoon} in $1964$ as an analogue of the extremal problem. In $1986$, K\'aszonyi and Tuza~\cite{kaszonyituza} characterized saturation functions for graphs by showing that a graph $G$ has a bounded saturation function if and only if $G$ has an isolated edge, otherwise, its saturation function is linear. Since then, saturation problems have become an important topic in extremal combinatorics. For an updated survey, see~\cite{graphsatsurvey}. 

In recent years, significant progress has been made in studying the saturation for $0$-$1$ matrices, which are equivalent to vertex-ordered bipartite graphs. Brualdi and Cao~\cite{brualdicao} were the first to investigate their saturation functions. Shortly after that, Fulek and Keszegh~\cite{fulkesz} showed that the dichotomy holds for $0$-$1$ matrices as well, identifying large families of matrices with linear saturation functions. On the other hand, Geneson~\cite{geneson} found an infinite family of matrices with a bounded saturation function. Finally, Berendsohn~\cite{berendsohn} gave a complete characterization of permutation matrices with bounded saturation functions. 

Inspired by these interesting results about $0$-$1$ matrices, the authors~\cite{boskesz023} initiated the study of saturation problems on vertex-ordered graphs. They considered two different orders on vertex sets, linear and cyclic order. They proved a dichotomy in both cases and  found infinite families of graphs that have bounded (respectively, linear) saturation functions. However, even for perfect matching graphs (which are the natural counterparts of permutation matrices from the $0$-$1$ matrix setting), the complete characterization is still unknown. 

On the other hand, saturation problems for directed graphs have not yet been studied systematically. F\"uredi et al.~\cite{Furedi1998MinimalOG} showed in a different context that $\overrightarrow{C_3}$, a directed cycle of length three, has a superlinear semisaturation function. More precisely, $ssat(n, \overrightarrow{C_3}) = n \log n + O(n \log \log n)$. This implies that the saturation function is also superlinear. In~\cite{graphsatsurvey} they ask for a matching upper bound for $sat(n, \overrightarrow{C_3})$. In fact, leaving out the two special vertices from the construction of~\cite{Furedi1998MinimalOG} shows that $sat(n, \overrightarrow{C_3}) = n \log n + O(n \log \log n)$ as well, answering this question.\footnote{This construction saturating $\overrightarrow{C_3}$ is the following. Let $k$ be an appropriate even number. We take a complete bipartite graph with parts of size $2k$ and at most $\binom{2k}{k}$ for an appropriate $k$. On the right side, the vertices correspond to different sets of size $k$ of the left side, and then we orient an edge from $w$ towards $v_E$ corresponding to some set $E$ if and only if $w \in E$. It is easy to check that this graph avoids $\overrightarrow{C_3}$, as even without orientation it does not contain a triangle. Also, adding any new edge with any orientation introduces a copy of $\overrightarrow{C_3}$.} This was observed by Alon and Fox~\cite{perscommalonfox}. 

Pikhurko~\cite{Pikhurko} considered a saturation problem on the class of cycle-free directed graphs, i.e. adding any missing edge introduces a cycle or a copy of a directed graph $G$, and showed that any family $\mathcal{F}$ of cycle-free directed graphs has $sat(n, \mathcal{F}) = O(n)$ with respect to this definition.

The study of extremal problems of edge-ordered graphs also leads to interesting problems. Their systematic study was started by Gerbner et al.~\cite{gerbner}. They proved an Erd\H os-Stone-Simonovits type theorem, which gives the exact asymptotic (a quadratic function) for the extremal function if the so-called order chromatic number of the forbidden graph is at least $3$. Thus, graphs with an order chromatic number $2$ are the most interesting to study further. They systematically took account of edge-ordered paths with $4$ edges.

This study was continued by Kucheriya and Tardos~\cite{kucheriyatardos1, kucheriyatardos2}. They characterized connected edge-ordered graphs that have a linear extremal function and showed that for every other connected edge-ordered graph the extremal function is $\Omega(n\log n)$, a dichotomy similar to that for vertex-ordered graphs. They also showed that for edge-ordered forests of order chromatic number $2$, the extremal function is $n\cdot 2^{O(\sqrt{{\log n}})}$. Note that if an edge-ordered graph contains a cycle, then its extremal function is $\Omega(n^c)$ for some $c>1$.

The motivation to study extremal problems of vertex-ordered graphs came from combinatorial geometry, where several natural problems can be described in such a way. It turns out that edge-ordered graphs have similar applications, as shown already in~\cite{gerbner}, and later in~\cite{xmon,curves}.

Continuing this line of research, we hereby initiate the study of saturation problems of edge-ordered graphs.

\subsection{Our results}

For (unordered) graphs~\cite{kaszonyituza}, $0$-$1$ matrices~\cite{fulkesz}, and vertex-ordered graphs~\cite{boskesz023}, it was shown that the saturation functions are either $O(1)$ or $\Theta(n)$. 	
Our main and somewhat surprising result is that there is no such dichotomy for edge-ordered graphs, similarly to directed graphs. We first consider the definition $sat_m$, and then the remaining two.

We identify infinitely many edge-ordered graphs with a superlinear saturation function. More precisely, we showed that for these graphs, the saturation function lies between $\Omega (n \sqrt{\log n})$ and $O\left(n \frac{\log n}{\log \log n}\right)$. For further families of graphs, we were able to show a weaker upper bound $O(n \log n)$, yet we were not able to find a non-trivial (subquadratic) upper bound that holds for every graph.

We also find infinite families of graphs with bounded and linear saturation functions. Among others, if $T$ is an edge-ordered tree, then $sat_m$ of $e_0 + T$ is bounded, where $e_0+T$ denotes the graph we get from $T$ by adding an isolated edge with the minimal label.

Dichotomy holds for the $sat_m$ definition in the following weak sense. If the minimal edge $e_0$ is isolated, then $sat_m$ is either bounded or linear. It is even possible, but we could not prove it, that the latter case never happens. Otherwise, if $e_0$ is not isolated, then $sat_m$ is at least linear. A summary of our results for the $sat_m$ definition is shown in Table \ref{table:sat}.

Next, we concentrate on the semisaturation problem of edge-ordered graphs. Unlike $0$-$1$ matrices~\cite{fulkesz} and vertex-ordered graphs~\cite{boskesz023}, where it was relatively simple to fully characterize semisaturation functions, for edge-ordered graphs, it appears to be more delicate. 
Our main result on semisaturation functions is a general upper bound $O(n \log n)$. We are also able to characterize the edge-ordered graphs with bounded $ssat_m$, they are the ones with isolated $e_0$. 
When $e_0$ is not isolated, then among others, we show that if $e_0$ is not participating in a triangle, then the graph has linear $ssat_m$.
A summary of our semisaturation results for $ssat_m$ definition is shown in Table \ref{table:ssat}.

We then consider the other two proposed variants: $sat_e$ and $sat_s$. We show several results following the proofs for the $sat_m$ definition. First, we show a weak dichotomy for $sat_e$, then characterize graphs with a bounded semisaturation function $ssat_e$. Moreover, we notice that the general semisaturation upper bound also holds for this definition. Later on, we find infinite families with bounded, linear, and superlinear (semi)saturation functions. 

We study the family of edge-ordered graphs $e_0 + G_k + e_{\max}$, where $e_0$ is the minimal edge and $e_{\max}$ is the maximal edge, and the underlying graph of $G_k$ is a complete graph $K_k$. It turns out that for some edge-orderings of $G_k$, the saturation function $sat_e$ is superlinear. This result shows that, even though they are related, the functions $sat_e$ and $sat_m$ behave differently, since for $sat_m$ we show that $sat_m(n, e_0 + G) = O(n)$ holds for any edge-ordered graph $G$. 
The summary of the results regarding $sat_e$ and $ssat_e$ functions can be seen in Table \ref{table:sate}.

Considering $sat_s$, by definition, most of the results about $sat_s$ are more general than those for $sat_m$, because any condition imposed on the minimal edge $e_0$ can be replaced by an equivalent condition on an arbitrary edge in a given graph. Even though we do not study this variant in detail, studying the $sat_s$ of an edge-ordered complete graph appears to be useful in discovering the difference between $sat_m$ and $sat_e$ saturation functions.

\subsection{Tables of results}

\begin{table}[ht]
	\centering
	\begin{tabular}{ | p{6.7cm} | p{3.6cm}| p{2.9cm} | } 
		\hline
		Graph & $sat_m$ & Reference \\
		\hline
		$e_0+G$  & $O(1)$ or $\Theta(n)$  & Theorem \ref{isolated_edge} \\
		\hline
		$e_0$ is not isolated  & $\Omega(n)$  & Theorem \ref{isolated_edge} \\
		\hline
		$G$, $N_G[a]=N_G[b]$, $G - e_0$ is bipartite  & $O(n \log n)$  & Theorem \ref{sat nlogn} \\
		\hline
		$v,w \in N_G(a) \cap N_G(b)$, $l(av) < l(aw)$ and $ l(bv) > l(bw)$  & $\Omega(n\sqrt{\log n})$  & Theorem \ref{D_0 general} \\
		\hline
		$D_0, D_1, D_2$  & $\Omega(n\sqrt{\log n})$  & Corollary \ref{D_0 lower bound} \\
		\hline
		$D_0+G$, $l(D_0) < l(G)$  & $O\left(n \frac{\log n}{\log \log n}\right)$  & Corollary \ref{D_0+G upper bound} \\
		\hline
		$D_3, D_4, D_5$  & $\Theta(n)$  & Claim \ref{D3-D5} \\
		\hline
		$C_k$, $k \geq 5$  & $\Theta(n)$  & Claim \ref{eocycles} \\
		\hline
		$e_0 + K_r, r \geq 2$  & $O(1)$  & Corollary \ref{e_0 complete} \\
		\hline
		$e_0 + F$, for monotone forest $F$  & $O(1)$  & Corollary \ref{mon forests} \\
		\hline
	\end{tabular}
	\caption{Review of saturation results, for missing definitions see the references.}
	\label{table:sat}
\end{table}

\begin{table}[ht]
	\centering
	\begin{tabular}{ | p{5.8cm} | p{4.5cm}| p{2.9cm} | } 
		\hline
		Graph & $ssat_m$ & Reference \\
		\hline
		$e_0 + G$  & $O(1)$  & Claim \ref{iso-edge-ssat} \\
		\hline
		$e_0$ is not isolated & $\Omega(n)$  & Claim \ref{iso-edge-ssat} \\
		\hline
		any $G$  & $O(n \log n)$  & Theorem \ref{nlogn upper bound} \\
		\hline
		$ v,w \in N_G(a) \cap N_G(b)$, $l(av) < l(aw)$ and $ l(bv) > l(bw)$  & $\Omega(n\sqrt{\log n})$  & Theorem \ref{D_0 general} \\
		\hline
		$D_0, D_1, D_2$  & $\Omega(n\sqrt{\log n})$  & Corollary \ref{D_0 lower bound} \\
		\hline
		$D_0+G$, $l(D_0) < l(G)$  & $O\left(n \frac{\log n}{\log \log n}\right)$  & Corollary \ref{D_0+G upper bound} \\
		\hline
		$e_0 = ab$, $N_G(a) \cap N_G(b) = \emptyset$  & $\Theta(n)$  & Claim \ref{disjoint_neighborhood} \\
		\hline
	\end{tabular}
	\caption{Review of semisaturation results.}
	\label{table:ssat}
\end{table}

\begin{table}[ht]
	\centering 
	\begin{tabular}{ | p{4.4cm} | p{6.6cm}| p{2.2cm} | } 
		\hline
		$e_0$ and $e_{\max}$ isolated  & $sat_e(n,G) = O(1)$ or $\Omega(n)$  & Theorem \ref{e_0e_max-isolated} \\
		\hline
		$e_0$ and $e_{\max}$ isolated & $ssat_e(n,G) = O(1)$ & Theorem \ref{bounded-ssat_e} \\
		\hline
		$e_0$ or $e_{\max}$ not isolated & $ssat_e(n,G) = \Omega(n)$  & Theorem \ref{bounded-ssat_e} \\
		\hline
		any $G$  & $ssat_e(n,G) = O(n \log n)$  & Corollary \ref{nlogn upper bound ssat_e} \\
		\hline
		$ab\in\{e_0,e_{\max}\}$,
		$v,w \in N_G(a) \cap N_G(b)$; $l(av) < l(aw)$ and $ l(bv) > l(bw)$  & $ssat_e(n,G) = \Omega(n\sqrt{\log n})$  & Corollary \ref{D_0 general sat_e} \\
		\hline
		$D_0$ & $sat_e(n, D_0) = O\left(n \frac{\log n}{\log \log n}\right)$   & Theorem \ref{D_0 sate upper bound} \\
		\hline
		$G_k$ certain edge-ordering of $K_k$, $k \ge 5$  & $ssat_e(n,e_0 + G_k + e_{\max}) = \Omega (n\sqrt{\log n})$  & Theorem \ref{complete_superlinear} \\
		\hline
		$C_{2k+1}, k \geq 1$  & $sat_e(n,C_{2k+1}) = \Theta(n) \newline ssat_e(n,C_{2k+1}) = \Theta(n)$  & Corollary \ref{eeocycles} \\
		\hline
		$M_k$ matching graph  & $sat_e(n,M_k) = O(1)$  & Corollary \ref{matchings} \\
		\hline
	\end{tabular}
	\caption{Results for $sat_e$ and $ssat_e$ definitions.}
	\label{table:sate}
\end{table}

\section{Dichotomies and general upper bounds}

\subsection{\texorpdfstring{$sat_m$}{sat?}}

In the vertex-ordered setting, it was possible to show that the saturation function of any graph is either bounded or linear~\cite{boskesz023}. Such a dichotomy does not hold for edge-ordered graphs in general, as we will see later. However, we can still separate between bounded and at least linear saturation functions. In addition, if the minimal edge is isolated, then we do have a dichotomy. Note that the isolated edge case was the only interesting case for vertex-ordered graphs~\cite{boskesz023}, as the remaining vertex-ordered graphs all have linear saturation function.

\begin{thm} \label{isolated_edge}
	Let $G$ be an edge-ordered graph. If the minimal edge $e_0$ of $G$ is isolated, then either $sat_m(n,G) = O(1)$ or $sat_m(n,G) = \Theta(n)$.
	Otherwise, if $e_0$ is not isolated, then $sat_m(n,G) = \Omega(n)$.
\end{thm}
\begin{proof}
	Assume first that $e_0$ is not isolated. There cannot be two isolated vertices in a saturating host graph $H$, otherwise by connecting them we will not introduce a new copy of $G$. Thus, $sat_m(n,G) \ge n/2-1$ and the result follows.
	
	Now, we assume that $e_0$ is isolated. To prove that $sat_m(n,G) = O(n)$ for any $G$ it is enough to take $H$ to be the union of $G - e_0$ and $n-k+2$ isolated vertices, where $k$ is the number of vertices of $G$. If we add an edge with a minimal label between any two isolated vertices, a copy of $G$ will be obtained. Next, we can create a saturating host graph $H'$ by greedily adding edges to $H$ (a common trick for other structures as well, e.g., for vertex-ordered graphs), such that $H'$ can have at most $\binom{n}{2} - \binom{n-k+2}{2} \le kn$ edges. This implies that $sat_m(n,G) = O(n)$. 
	
	On the other hand, if $G$ has a sublinear saturation number, then for every large enough $n_0$ there is a host graph $H_0$ on $n_0$ vertices that has two isolated vertices. Then it is easy to see that by adding isolated vertices to $H_0$, to ensure it has $n$ vertices in total, we still must have a saturating host graph. Therefore, we conclude that if the saturation number is not linear, then it has to be bounded.
\end{proof}

It is actually possible that if $e_0$ is isolated, then $sat_m$ is always bounded. We can prove this in the case of semisaturation.

\begin{claim} \label{iso-edge-ssat}
	Let $G$ be an edge-ordered graph. If the minimal edge $e_0$ of $G$ is isolated, then $ssat_m(n,G) = O(1)$, otherwise $ssat_m(n,G) = \Omega (n)$. 
\end{claim}
\begin{proof}
	Assume that $e_0$ is not isolated. If there are two isolated vertices in the host graph $H$, by connecting them we do not introduce a new copy of $G$. Thus, $sat_m(n,G) \ge n/2-1$. This is the same proof as in Theorem \ref{isolated_edge}.
	
	Assume now that $e_0$ is isolated. Take $H$ to be a graph on $n$ vertices that has two disjoint copies of $G$ and all the remaining vertices isolated. Then we add edges between non-isolated vertices to get a host graph with a bounded number of edges.
\end{proof}

We now proceed to prove a general semisaturation upper bound. We denote by $N_G(a)$ (resp. $N_G[a]$) the open (resp. closed) neighborhood of $a$ in $G$.

\begin{thm} \label{nlogn upper bound}
	For every edge-ordered graph $G$, we have $ssat_m(n,G) = O(n \log n)$.
\end{thm}

\begin{proof}
	Let $G$ be an edge-ordered graph with the corresponding labeling function $l: E(G) \rightarrow \{0, 1, ..., m-1\}$ and denote by $e_0 = ab$ its minimal edge, that is, $l(e_0) = 0$. We define an additional labeling function for host graphs $L: E(H) \rightarrow \mathbb{N}^3$, where $\mathbb{N}^3$ is equipped with lexicographic order. We write $x.y.z$ for $(x,y,z) \in \mathbb{N}^3$, and $x.y$ for $(x,y,0)$. Let $N_G(a) \setminus b := \{a_1, \dots, a_r\}$ and $N_G(b) \setminus a := \{b_1, \dots, b_s\}$. Let  $k = \lceil \log n \rceil$ and $n'< n$ to be determined later. Consider $k$ disjoint copies of the graph $G^{ab} = G \setminus \{a,b\}$, which we denote by $G^{ab}_i$, for $1 \leq i \leq k$. For every edge $e \in E(G^{ab})$, we label the corresponding edge $e_i \in E(G^{ab}_i)$ by: \[L(e_i) = l(e).i, \text{ for all } 1 \leq i \leq k.\] Furthermore, we add $n'$ isolated vertices $V = \{v_1, ..., v_{n'}\}$. 
	
	Take $k$ bipartitions $(X_i,Y_i)$ of $V$ such that every pair of vertices is separated in at least one of them. This can be done because $n'\le 2^k$. For example, the $i$-th bipartition is obtained according to the $i$-th digit in the binary form of the vertex indices. If necessary, some bipartitions can be repeated.
	
	Fix some $i\le k$. Then, take all the vertices $a_{i,j}$ in $G^{ab}_i$ that correspond to the neighbors of $a$. We add an edge between $v_t \in X_i$ and $a_{i,j}$, where $j \in \{1,...,r\}$. We label such edges with the following function:
	\[L(v_ta_{i,j}) = l(aa_j).i.t, \text{ for all } 1 \leq j \leq r.\]
	
	Analogously, we do the same for the neighbors of $b$, but this time by connecting them with isolated vertices $v_t$ in $Y_i$. We label them by:
	\[ L(v_tb_{i,j}) = l(bb_j).i.t, \text{ for all } 1 \leq j \leq s.\]
	
	This graph has $O(kn')$ edges. We can see that by adding an edge between $x_i \in X_i$ and $y_i \in Y_i$, the graph induced by $G_i^{ab}$ and these two vertices is isomorphic to $G$, where the edge $x_iy_i$ plays the role of the minimal edge $ab$. 
	
	As every pair of vertices in $V$ is separated in some bipartition, adding an edge between any two vertices in $V$ creates a copy of $G$. Set $n'$ such that this graph has exactly $n$ vertices. By greedily adding edges between vertices that are not in $V$, we obtain a host graph that semisaturates $G$ and has $O(n \log n)$ edges.	
\end{proof}

While the upper bound $O(n \log n)$ could also be true for $sat_m$, we can prove it only for a subclass of graphs.

\begin{thm} \label{sat nlogn}
	Let $G$ be an edge-ordered graph such that $G-e_0$ is bipartite and $N_G[a] = N_G[b]$, where $e_0 = ab$ is the minimal edge. Then \[sat_m(n,G) = O(n \log n).\]
\end{thm}

\begin{proof}
	We construct the graph $H$ the same way as in Theorem \ref{nlogn upper bound} and we show that $H$ is $G$-free. First, observe that if $N_G[a] = N_G[b] = \{a,b\}$, then $e_0$ is an isolated edge and the result follows by Claim \ref{isolated_edge}. Otherwise, there exists a vertex $w$ that creates a triangle with $a$ and $b$, which means that $G$ is non-bipartite. We will show that $H$ given in Theorem \ref{nlogn upper bound} is bipartite and it follows that $H$ avoids $G$. Then we again greedily extend $H$ to get a saturating host graph for $G$.
	
	Let $(X, Y)$ be a bipartition of $G-e_0$ (which exists by the assumption of the theorem), then notice that $a$ and $b$ should be in the same part, say $X$, while all of their neighbors $N^{ab} := N_G(a) \cap N_G(b)$ have to be in $Y$.
	Since $G-e_0$ is bipartite, then $G_i^{ab}$ is also bipartite for every $1 \leq i \leq k$. We now can construct a bipartition $(X_H, Y_H)$ of $H$, and therefore prove the claim. We take a disjoint union of all $G_i^{ab}$, such that the corresponding neighbors $N_i^{ab}$ of $a$ and $b$ are all in $Y_H$. Then we add the set $V$ as defined in the proof of Theorem \ref{nlogn upper bound}, to the set $X_H$ and we connect all the vertices in $V$ with all the vertices in $\bigcup_{i=1}^k N_i^{ab}$, which gives a desired bipartition of $H$. So, we can conclude that $H$ is indeed $G$-free.
\end{proof}

\subsection{\texorpdfstring{$sat_e$}{sat?}}

We are only able to prove a weak dichotomy for $sat_e$. We denote by $e_0$ the minimal edge and by $e_{\max}$ the maximal edge of $G$.

\begin{thm} \label{e_0e_max-isolated}
	Let $G$ be an edge-ordered graph. If $e_0$ and $e_{\max}$ are isolated, then $sat_e(n,G) \in \{ O(1), \Omega(n)\}$. Otherwise, $sat_e(n,G) = \Omega(n)$.
\end{thm}
\begin{proof}
	We follow the same ideas as in the proof of Theorem \ref{isolated_edge}. First, we consider the case when either $e_0$ or $e_{\max}$ is not isolated. Without loss of generality, $e_0$ is not isolated. Then adding a missing edge with the minimal label between two isolated vertices of $H$ does not create a copy of $G$. Therefore, $H$ cannot contain two isolated vertices, which implies $sat_e(n,G) \geq n/2-1$. 
	
	Let $e_0$ and $e_{\max}$ be isolated, and assume that $G$ has a sublinear saturation function. For a large enough $n$, the host graph must have at least two isolated vertices. Since they can be replaced by arbitrarily many isolated vertices without affecting the saturation property, we conclude that if the saturation function is sublinear, it must in fact be bounded.		
\end{proof}

Notice that in the $e_0$ and $e_{\max}$ isolated case when $sat_e$ is unbounded, unlike in Theorem \ref{isolated_edge} where we proved $\Theta(n)$, here we only have $\Omega(n)$. The reason for this is detailed later.

Similarly to Claim \ref{iso-edge-ssat}, for semisaturation we can prove a stronger statement.

\begin{thm} \label{bounded-ssat_e}
	Let $G$ be an edge-ordered graph. Then \[ssat_e(n,G) = \begin{cases}
		O(1) & \text{if $e_0$ and $e_{\max}$ are isolated,} \\
		\Omega (n) & \text{otherwise.} \\ \end{cases}\]
\end{thm}
\begin{proof}
	Let $H$ be a host graph on $n$ vertices and $m$ edges, such that $H = e_1 + G_1 + G_2 + G_3 + G_4 + e_m + V$, where $G_i$ are copies of $G - \{e_0, e_{\max}\}$, so that $l(e_1) < l(G_1) < l(G_2) < l(G_3) < l(G_4) < l(e_m)$ and $V$, the set of the remaining vertices of $H$, is a set of isolated vertices.  Let $e'=uv$ be an edge that we can add to $H$ with $u\in V$ isolated. We distinguish two cases:
	\begin{enumerate}
		\item $l(e') > l(G_2)$. If $v \in V(G_2 + G_3 + G_4 + e_m + V)$, then $e_1 + G_1 + e'$ is a copy of $G$. Otherwise, if $v \in V(e_1 + G_1)$, then  $e'' + G_2 + e'$ is a copy of $G$, where $e'' \in E(e_1 + G_1)$, unless $G_1$ contains no edges. In that case, the graph $e' + G_2 + e_m$ is a copy of $G$.
		\item $l(e') < l(G_3)$. Similarly, if $v \in V(e_1 + G_1 + G_2 + G_3 + V)$, then $e' + G_4 + e_m$ is a copy of $G$. Otherwise, if $v \in V(G_4 + e_m)$, then $e' + G_3 + e''$ is a copy of $G$, where $e'' \in E(G_4+e_m)$, unless $G_4$ contains no edges. In that case, the graph $e_1 + G_3 + e'$ is a copy of $G$.
	\end{enumerate}
	
	Thus, adding edges greedily to $H$ that are not in $V$, we get a host graph that semisaturates $G$ with $O(1)$ edges.
	
	Otherwise, if for example, $e_0$ is not isolated, then by adding any new edge with minimal label between two vertices in $V$, we do not create a copy of $G$. Therefore, the semisaturation function is at least linear.
\end{proof}

We are also able to show the general upper bound for semisaturation. Let $G^{\text{rev}}$ be the edge-ordered graph obtained from $G$ such that the underlying graph remains the same, but the linear order of its edges is reversed. By Theorem \ref{nlogn upper bound} we get that $sat_m(n,G) = O(n \log n)$ and $sat_m(n,G^{\text{rev}}) = O(n \log n)$. By gluing the host graphs $H$ and $H^{\text{rev}}$ on the same independent set $V$ such that $l(H^{\text{rev}}) < l(H)$, we get a new host graph which semisaturates $G$ in $sat_e$ sense for any missing edge in $V$.

\begin{cor} \label{nlogn upper bound ssat_e}
	For every edge-ordered graph $G$, we have $ssat_e(n,G) = O(n \log n)$.
\end{cor}

Unlike $sat_m$ definition, we do not know whether it is possible to show this upper bound for the saturation function of some family similar to the one described in Theorem \ref{sat nlogn}.

\subsection{\texorpdfstring{$sat_s$}{sat?}}

Recall that $sat_s(n, G) \leq sat_m(n, G)$ for every edge-ordered graph $G$. So, all the results about upper bounds can be applied to the $sat_s$ definition.
The forthcoming results are more general than those for the $sat_m$ definition, which is related to the fact that we can usually replace the assumptions about the minimal edge $e_0$ with any edge in $G$. We state here the weak dichotomies and the general semisaturation upper bound. We do not include the proofs of these results, instead, we provide references for the corresponding $sat_m$ results, whose proofs can be adapted to this framework.
In Section \ref{sats section}, we will return to this definition and explain several other results.

The first result is the analogue of Theorem \ref{isolated_edge}. Note that it does not require an edge with a minimal or maximal label to be isolated.

\begin{thm} \label{weak dichotomy sats}
	Let $G$ be an edge-ordered graph. If $G$ contains an isolated edge, then either $sat_s(n,G) = O(1)$ or $sat_s(n,G) = \Theta(n)$. Otherwise, $sat_s(n,G) = \Omega(n)$.
\end{thm}

We continue with the weak dichotomy of semisaturation functions. The proof is completely analogous to the proof of Claim \ref{iso-edge-ssat}.

\begin{claim} \label{bounded-ssat_s}
	Let $G$ be an edge-ordered graph. Then,
	\[ssat_s(n,G) = \begin{cases}
		O(1) & \text{if $G$ contains an isolated edge,} \\
		\Omega (n) & \text{otherwise.} \\
	\end{cases}\]
\end{claim}

The general upper bound for semisaturation follows directly from Theorem \ref{nlogn upper bound}.

\begin{cor} \label{nlogn upper bound ssats}
	For every edge-ordered graph $G$, we have $ssat_s(n,G) = O(n \log n)$.
\end{cor}

Many of our results about the $sat_m$ function can be replicated for the $sat_s$ function by possibly changing the assumptions. For example, the following result follows by a simple adaptation of the proof of Theorem \ref{sat nlogn}.

\begin{claim} \label{sats nlogn}
	Let $G$ be an edge-ordered graph. If there exists an edge $ab \in E(G)$ such that $G-ab$ is bipartite and $N_G[a] = N_G[b]$. Then $sat_s(n,G) = O(n \log n)$.
\end{claim}

\section{Unbounded saturation functions}

In this section, we deal only with $sat_m$ and $ssat_m$, while the other variants are discussed in later sections.

\subsection{Bipartite coverings}

In order to prove the existence of edge-ordered graphs with superlinear saturation functions, we will need an auxiliary result about bipartite coverings of almost complete graphs. 

A \emph{bipartite covering} of a graph $G$ is a collection of bipartite graphs such that each edge of $G$ belongs to at least one of them. The \emph{capacity} of the covering is the sum of the numbers of vertices of these bipartite graphs.	

\begin{thm}\label{thm:bip}\cite{katona}
	Let G be a graph on $n$ vertices and let $d_1, d_2, . . . , d_n$ denote the degrees of its vertices. Then the capacity of any bipartite covering of $G$ is at least
	\[\sum_{i=1}^{n}(\log n- \log (n-d_i)).\]
\end{thm}

This was shown by Katona and Szemer\'edi~\cite{katona}, see also the paper of Alon~\cite{alon} for a generalization which also has a similar and short proof. See also the exact result of Dong and Liu~\cite{dong007} that in a bipartite covering of the complete graph $K_n$ there is a vertex which has incident edges from at least $\log n$ of the bipartite graphs.

\begin{cor}\label{cor:bipdense}
	Let $G$ be a graph on $n$ vertices \footnote{As usual, $n$ is assumed to be large enough.} with at most $cn^{2-\epsilon}$ missing edges for some constants $c,\epsilon > 0$. Let $\cal B$ be a bipartite covering of $G$, then there is a vertex which has incident edges from at least $c' \log n$ of the bipartite graphs in $\cal B$, where $c'$ depends only on $c$ and $\epsilon$.		
\end{cor}

\begin{proof}
	Let $x=2cn^{1-\epsilon}$. Remove from $G$ all vertices whose co-degree (i.e., their degree in the complement graph) is at least $x$, deleting at most $n/2$ vertices in total. Let $G'$ be the resulting graph on $n' \ge n/2$ vertices. Moreover, every vertex in $G'$ has co-degree less than $x$, which implies that its degree is at least $n'-x-1$. 
	
	By Theorem \ref{thm:bip}, the capacity of $\cal B$ restricted to $G'$ is at least \[\sum_{i=1}^{n'}(\log n'- \log (n'-d_i))> \frac{n}{2}(\log n'- \log (x+1))= \frac{n}{2}(\epsilon\log n-O(1)),\]
	
	where the $O$ notation hides dependence on $c$ and $\epsilon$. Therefore, one of the vertices must be in at least $c'\log n$ bipartite graphs for some constant $c'$.
\end{proof}

\subsection{Superlinear saturation functions}

In this subsection, we show a family of graphs that have superlinear saturation and semisaturation functions, showing that there is no dichotomy. 
Notice that, by Theorem \ref{D_0 general} dichotomy does not hold either for saturation or for semisaturation. This most likely suggests that studying semisaturation of edge-ordered graphs might be more demanding than it was for vertex-ordered graphs~\cite{boskesz023}, where it was possible to fully  characterize semisaturation functions.

\begin{thm} \label{D_0 general}
	Let $G$ be an edge-ordered graph such that its minimal edge is $ab$. If there exist two vertices $v,w \in N_G(a) \cap N_G(b)$ such that $l(av) < l(aw)$ and $ l(bv) > l(bw)$, then \[sat_m(n,G) \geq ssat_m(n,G)=\Omega(n\sqrt{\log n}).\]
\end{thm}

\begin{figure}[ht]
	\centering
	\includegraphics[width=0.7\textwidth]{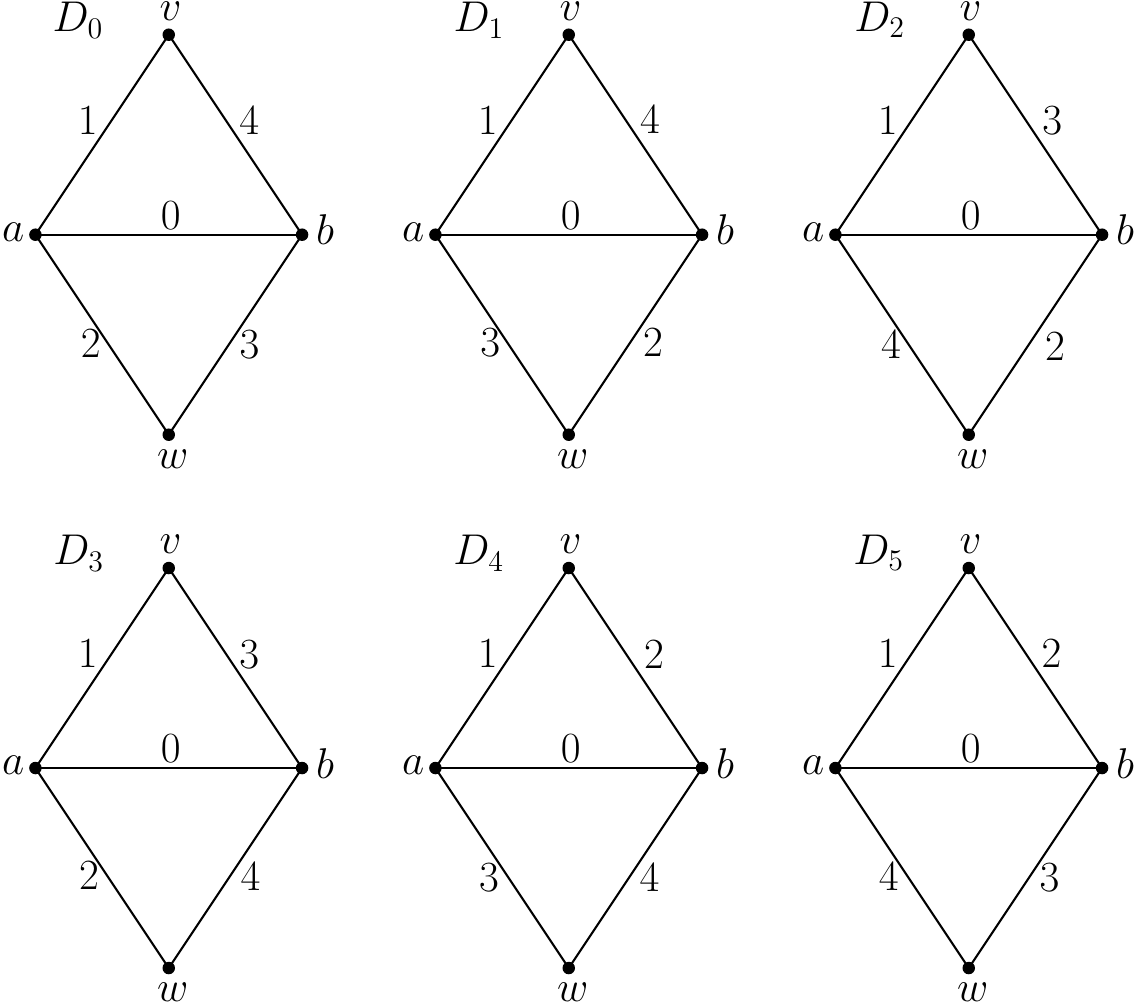}
	\caption{\label{fig:D012345} Six edge-ordered diamond graphs whose minimal edge is $ab$.}
\end{figure}

\begin{proof}
	We note that by the assumptions, the minimal edge $ab$ is not isolated.
	Let $H$ be a saturation host graph on $n$ vertices. For a pair of vertices that do not span an edge of $H$, denoted by $f$, adding an edge at $f$ with the minimal label must create a copy of $G$. Thus there exists (at least one) pair of vertices $e$ such that the four vertices $f\cup e$ span a copy of the structure on four vertices and five edges which is required by our assumption in which $f$ plays the role of $ab$. Choosing one such $e$ for each $f$, this defines an assignment function $\pi(f)=e$ on the pairs of vertices of $H$ that do not span an edge. 
	
	Because of how $f$ and $e$ correspond to vertices in $G$, for one vertex of $f$ we have that the labels of the two edges going from this vertex to the two vertices of $e$ are in reverse order compared to the labels of the two edges going from the other vertex of $f$ to the vertices of $e$. 
	See, e.g., in the case $G=D_0$ in Figure \ref{fig:D012345}, where $f$ corresponds to $ab$ and $e$ corresponds to $vw$. Therefore, for every $e$ we can split the vertices connected to both vertices of $e$ into $2$ types according to the order of labels and $e$ can be assigned only to pairs $f$ which contain one vertex from each type. That is, for every pair of vertices $e$ there is a bipartition of the vertices such that all non-edges in $\pi^{-1}(e)$ lie between the two parts.
	
	Assume on the contrary that $H$ has less than $cn\sqrt{\log n}$ edges for some $c\le 1$ to be chosen later. The non-edges of $H$ must be covered by the sets of non-edges $\pi^{-1}(e)$ where $e$ is a pair of vertices in $H$. As each $\pi^{-1}(e)$ is bipartite, we apply Corollary \ref{cor:bipdense} to the complement of $H$. First, we choose constants $c_0,\epsilon_0$ independently from $c$ such that $ n\sqrt{\log n}\le c_0(n/2)^{2-\epsilon_0}$ holds. 
	
	Applying Corollary \ref{cor:bipdense} with these constants, we find a vertex $v$ participating in at least $c'\log n$ bipartitions covering the non-edges of $H$. For each $e$ such that $\pi^{-1}(e)$ contains a non-edge that covers $v$, there is an edge in $H$ from $v$ to both vertices of $e$. Hence, $v$ is connected to at least $c'\log n$ such distinct pairs of vertices, therefore, to at least $\sqrt{c'\log n}$ vertices in $H$. Delete $v$ from $H$ and repeat this procedure $n/2$ times, at each step we apply Corollary \ref{cor:bipdense} with the same $c_0$ and $\epsilon_0$ to find and delete a vertex with degree at least $\sqrt{c'\log (n/2)}$. Altogether, we have found at least $\frac{n}{2} \sqrt{c'\log (n/2)}$ edges in $H$. When $c$ is small enough this contradicts our assumption that $H$ has less than $cn\sqrt{\log n}$ edges.
\end{proof}

\begin{cor} \label{D_0 lower bound}
	Let $D_0, D_1, D_2$ be the graphs shown in Figure \ref{fig:D012345}. Then \[ssat_m(n,D_i)=\Omega(n\sqrt{\log n}),\] for $i \in \{0,1,2\}$. 
\end{cor}

We would like to better understand which edge-ordered graphs are not covered by Theorem \ref{D_0 general}. We will discuss later the graphs where the end vertices of the minimal edge have disjoint neighborhoods. For now, we assume that $N_G(a) \cap N_G(b)$ is non-empty.

\begin{remark} \label{increasing_order}
	Let $\{u_1, u_2, ..., u_k\} := N_G(a) \cap N_G(b)$, for $k \geq 1$ and $ab$ is the minimal edge of $G$. We can assume wlog. that $l(au_1) < ... < l(au_k)$. If $l(bu_1) < ... < l(bu_k)$ does not hold, then there exist $1 \leq i < j \leq k$ such that $l(au_i) < l(au_j)$ and $l(bu_i) > l(bu_j)$, which by Theorem \ref{D_0 general} implies that $G$ has a superlinear saturation function.
\end{remark}

The following elementary result will be useful in the proof where we discuss the upper bound for the saturation function of the diamond graph $D_0$.

\begin{lem} \label{k!=n}
	Let $k, n \in \mathbb{N}$ such that $k! = n$. Then $k = O\left(\frac{\log n}{\log \log n}\right)$.
\end{lem}

\begin{proof}
	To describe $k$ in terms of $n$ we can use Stirling's formula to obtain $k \log k = \Theta(\log n)$, which is equivalent to
	\[\frac{c_1 \log n}{\log k} \leq k \leq \frac{c_2 \log n}{\log k},\]
	for some positive constants $c_1$ and $c_2$. From the first inequality we get that $c_1 \log n \leq k \log k \leq k^2$, which is equivalent to $\sqrt{c_1 \log n} \leq k$. If we plug this into the second inequality we get
	\[ k \leq \frac{c_2 \log n}{\log \sqrt{c_1 \log n}} = O\left(\frac{\log n}{\log \log n}\right). \]
\end{proof}

As an edge-ordered complete bipartite graph can encoded with its biadjacency matrix, in which every entry describes the corresponding label, we introduce the notion of \emph{sat property} that will be used throughout our next theorem.

\begin{defi}
	A $L = (\ell_{i,j})_{1 \leq i \leq k, 1 \leq j \leq n}$ satisfies the \emph{sat property} if the following holds:
	for every two column indices $1 \leq j_1 , j_2 \leq n$, there exist two row indices $1 \leq i_1 , i_2 \leq k$ such that we have $\ell_{i_2,j_1} < \ell_{i_1,j_1} < \ell_{i_1,j_2}< \ell_{i_2,j_2}$.
\end{defi}

\begin{figure}[ht]
	\centering
	\includegraphics[width=0.7\textwidth]{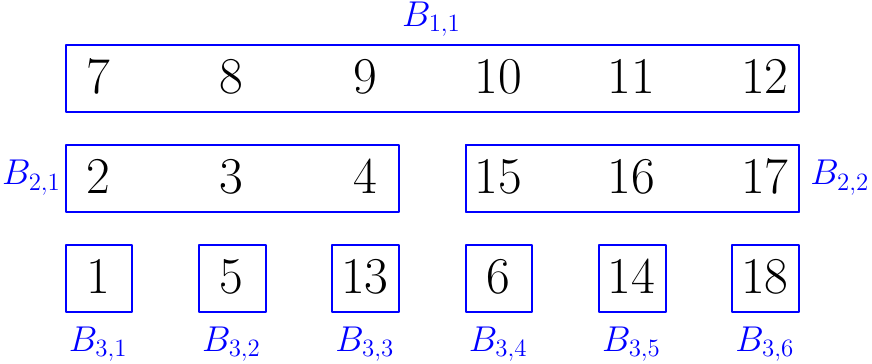}
	\caption{\label{fig:blocks} A matrix $L$ and its blocks for the case $n=6$.}		
\end{figure}

See Figure \ref{fig:blocks}, for an example of a matrix satisfying the sat property. We now proceed to establish an upper bound for $D_0$ which is close to its lower bound.

\begin{thm} \label{D_0 upper bound}
	\[ssat_m(n,D_0) \leq sat_m(n,D_0)=O\left(n \frac{\log n}{\log \log n}\right).\]
\end{thm}

\begin{proof}
	We show the upper bound by constructing a certain labeling of the complete bipartite graph $H = K_{k, n}$ (for appropriate small $k$) which will be a host graph for $D_0$. Since $D_0$ is non-bipartite, it follows that $H$ is $D_0$-free regardless of its labels. We will label the edges of $H$ such that it saturates $D_0$ for the non-edges within the $\bar K_n$ part, at the end we just greedily extend this graph with edges put within the $\bar K_k$ part to get a host graph for $D_0$. This time, we assume that the labeling function takes values in $\mathbb{R_{>}}$. 
	
	An alternative way to look at this problem is to describe the labeling in terms of the biadjacency matrix of $H$, which we denote by $L = (\ell_{i,j})_{1 \leq i \leq k, 1 \leq j \leq n}$, with entries from $\mathbb{R_{>}}$. We can check that by labeling the edges of $H$ according to a matrix $L$ that satisfies the sat property, we get a host graph that saturates $D_0$ for the non-edges within the $\bar K_n$ part of $H$ (the sat property is enough but not necessary for that). Notice that $L$ has the sat property if and only if for every pair of columns, the submatrix defined by these two columns has the sat property. See Figure \ref{fig:satproperty} for an illustration.
	
	\begin{figure}[ht]
		\centering
		\includegraphics[width=0.8\textwidth]{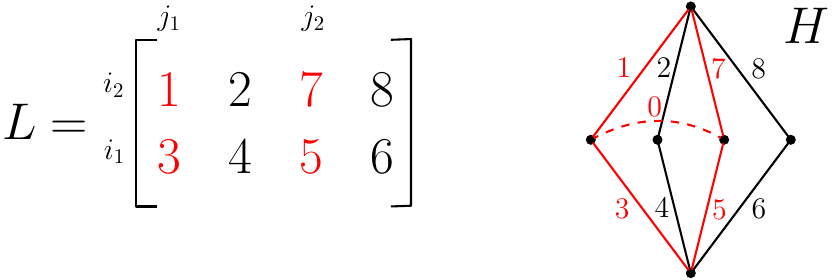}
		\caption{\label{fig:satproperty} The submatrix of $L$ defined by the columns $1$ and $3$ that satisfies the sat property. The host graph $H$ associated to $L$ that saturates $D_0$ for the corresponding dashed non-edge.}	
	\end{figure}
	
	First, let us assume that $n=k!$ for some $k \in \mathbb{N}$ and show that the matrix $L$ of size $k\times n$ which will be determined below satisfies the sat property. We partition the elements of the matrix $L$ into blocks $(B_{i,j})_{1 \leq i \leq k, \\ 1 \leq j \leq i!}$ such that the $i$-th row is equipartitioned into $i!$ blocks:
	\[B_{i,j} = \left\{\ell_{i,\frac{(j-1)n}{i!}+1}, \dots , \ell_{i,\frac{jn}{i!}}\right\}.\]
	
	Considering how the blocks in each row partition the columns, we can observe that this partition for any given row is a refinement compared to the row above, that is, for each block $B$ in the $(i-1)$th row, there are exactly $i$ blocks in the $i$-th row such that their columns equipartition the columns of $B$. See Figure \ref{fig:blocks} for an example.
	
	The labeling will be such that within any block the labels are in increasing order, that is, $\ell_{i,j_1} < \ell_{i,j_2}$ for any two entries in the same block whenever $j_1 < j_2$. Additionally, for any two blocks their labels can be separated, that is, for any two blocks $B_{i_1, j_1}, B_{i_2, j_2}$ all the elements of $B_{i_1, j_1}$ are smaller (or larger) than all the elements of $B_{i_2, j_2}$. We denote this by  $B_{i_1, j_1} < B_{i_2, j_2}$ (or $B_{i_1, j_1} > B_{i_2, j_2}$). 
	
	The labeling is defined row by row. First, we label the entries in $B_{1,1}$, the first row of the matrix, with any $n$ different positive real numbers in increasing order. Then we proceed to the second row by labeling $B_{2,1}$ and $B_{2,2}$ such that $B_{2,1} < B_{1,1} < B_{2,2}$. Note that if we take any two columns $j_1$ and $j_2$ such that $j_1 \leq \frac{n}{2} < j_2$, then these two columns and the first two rows of the matrix $L$ define a submatrix that has the sat property. We continue with the third row by labeling the blocks such that
	\[B_{3,1} < B_{2,1} < B_{3,2} < B_{1,1} < B_{3,3} \text{ and } B_{3,4} < B_{1,1} < B_{3,5} < B_{2,2} < B_{3,6}. \]
	See again Figure \ref{fig:blocks}, for an example when $n = 6$ and observe that it satisfies the sat property.
	
	In a general step, if we assume that $i-1$ rows are labeled for $i \geq 2$, then we label the $i$-th row as follows. We take any $j$ such that $1 \leq j \leq (i-1)!$ and we want to label all the blocks in row $i$ that share a column with $B_{i-1,j}$. These are exactly the blocks $B_{i,(j-1)i+1}, \dots, B_{i,ji}$ and we label them such that
	\[ B_{i,(j-1)i+1} < A_1 < B_{i,(j-1)i+2} < ... < A_{i-1} < B_{i,ji}, \]
	where $\{A_1, \dots, A_{i-1}\}$ are the blocks that share a column with $B_{i-1,j}$ (including $B_{i-1,j}$), (as there is exactly one such block in each row). Additionally, these are indexed in an increasing order, i.e., $A_1 < \dots < A_{i-1}$. In the example above, if we want to label the block $B_{3,2}$, then $A_1 := B_{2,1}, A_2 := B_{1,1}$.
	
	Next, we prove by induction on $i$ that if we take any two columns $j_1$ and $j_2$ that correspond to different blocks in the $i$-th row, then these two columns and the first $i$ rows of the matrix $L$ define a submatrix that has the sat property. The base case $i=2$ is clear. Suppose $i>2$ and the two columns correspond to different blocks of the row $i-1$, then this holds by induction. Finally, if $i>2$ and the two columns were separated first in the $i$-th column, say they are in blocks $B_1$ and $B_2$, where $B_1$ comes first (has smaller column indices). By the definition of labeling, we have $B_1<B_2$. Also, in each earlier row, there is a block that contains both columns, and there is one such block $A$ in some row $i'$ for which $B_1<A<B_2$. Since the entries within $A$ are increasing, the four entries in these two columns and in rows $i,i'$ show that the submatrix indeed has the sat property.
	
	When we reach the last row, we see that the blocks are all singletons, and we stop. As each column is separated in this last row, by the inductive claim, each pair of columns defines a submatrix with the sat property, and thus the whole matrix $L$ also has the sat property, as required.
	
	Thus, we obtain a labeling of $H$ from the matrix $L$. Finally, we greedily keep adding edges with minimal labels to $H$ within the $K_k$ part until the graph saturates $D_0$. As the number of edges in $H$ is at most $kn+\binom{k}{2}$. By Lemma \ref{k!=n}, we have $k = O\left(n \frac{\log n}{\log \log n}\right)$. This completes the proof under the assumption that $k! = n$.
	
	It remains to deal with the case when there is no $k \in \mathbb{N}$ such that $n=k!$. In this case set $k$ to be the smallest number with $n<k!$ and determine the $k\times n$ matrix $L$ the same way except that when defining the blocks in the $i$-th row, instead of equipartitioning the columns of each block in the $(i-1)$th row,  we now partition them to $i$ almost equal size blocks (i.e., the difference of their sizes is at most one). Note that since $(k-1)! < n$, then in $(k-1)$th row, each block is of size at least one. Only in the last row we may have blocks of size zero, but this does not change asymptotically the result and it is straightforward to check that the same proof works.
\end{proof}

Even though the bounds are not far, determining the exact order of magnitude of $sat_m(n,D_0)$ is left open. We can generalize this upper bound for further graphs. Before stating the result, we recall the notation: for two edge-ordered graphs $F_1$ and $F_2$, we write $l(F_1) < l(F_2)$ if all the labels of $F_1$ are smaller than all the labels of $F_2$.

\begin{lem}
	Let $F_1$ and $F_2$ be two edge-ordered graphs such that $sat_m(n, F_1) = O(f(n))$, where $f(n) = \Omega(n)$ and $l(F_1) < l(F_2)$. Then $sat_m(n, F_1 + F_2) = O(f(n))$.
\end{lem}
\begin{proof}
	Let $H(F_1)$ be the host graph for $F_1$, and assume $F_2$ has $k$ vertices. Then we take a union $H(F_1) + F_2^*$, where $F_2^*$ is a relabeled copy of $F_2$ such that $l(H(F_1)) < l(F_2^*)$. We argue that this construction indeed avoids $F_1+F_2$. Assume there is a copy of $F_1$ in $H(F_1) + F_2^*$, clearly it cannot be found in the part $H(F_1)$, so it has to have at least one edge in $F_2^*$. In that case $F_2$ has to be found in the rest of $F_2^*$, which is not possible because of its size. 
	
	Notice that by adding any edge within $H(F_1)$ we get a copy of $F_1+F_2$, and we can greedily add any other edge if it does not create a copy of $F_1+F_2$, which only adds at most $kn+k^2$ edges. Since $f(n)$ is at least linear, we get that $sat_m(n, F_1 + F_2) = O(f(n))$.
\end{proof}

Using this result, we can generalize the upper bound obtained for $D_0$ in Theorem \ref{D_0 upper bound}. Together with Theorem \ref{D_0 general}, we identify an infinite family of edge-ordered graphs with a superlinear saturation function.

\begin{cor} \label{D_0+G upper bound}
	Let $G$ be an edge-ordered graph, such that $l(D_0) < l(G)$. Then \[\Omega(n\sqrt{\log n}) = ssat_m(n,D_0+G)\le sat_m(n,D_0+G) = O\left(n \frac{\log n}{\log \log n}\right).\]
\end{cor}

\subsection{Linear saturation functions}

First, we discuss the remaining three graphs that have the diamond graph as their underlying graph.

\begin{claim} \label{D3-D5}
	For $i\in\{3,4,5\}$ we have $sat_m(n, D_i)= \Theta(n)$, where $D_i$ are the graphs as shown in Figure \ref{fig:D012345}.
\end{claim}

\begin{proof}
	By Theorem \ref{isolated_edge}, it is enough to show that $sat_m(n, D_i)= O(n)$ for all $i\in\{3,4,5\}$. Take $H$ to be a host graph on $n$ vertices such that its underlying graph is a complete bipartite graph $K_{2,n-2}$. Denote the vertices of $H$ by $u_1, u_2$ on one side and $v_1, ..., v_{n-2}$ on the other side. With a proper edge-ordering we will be able to define for each $D_i$ a saturating graph $H_i$ that has linearly many edges, from where the conclusion will follow. Notice also that since in all the cases the underlying graph of $H_i$ is bipartite, and $D_i$ is not, we can easily conclude that $H_i$ is $D_i$-free. Therefore, we will only need to check that $H_i$ actually saturates $D_i$.
	
	We start with $D_3$, the following labeling of $H$ gives $H_3$:
	\[ l(u_1v_i) = 2i-1 \text{ and } l(u_2v_i) = 2i \text{, for all } i \in \{1,2, ..., n-2\}. \]
	To check that $H_3$ saturates $D_3$, connect an edge $v_iv_j$ for any $1 \leq i < j \leq n-2$ and label it with zero. Then, notice that the graph induced by four vertices $u_1, u_2, v_i, v_j$ is isomorphic to $D_3$. We can add the edge $u_1u_2$ with the largest label, this does not introduce a copy of $D_3$ and does not affect the linearity of $sat_m(n, D_3)$.
	
	We continue with the case $D_4$, the following labeling of $H$ gives $H_4$:
	\[ l(u_1v_i) = i \text{ and } l(u_2v_i) = n-2+i \text{, for all } i \in \{1,2, ..., n-2\}. \]
	Again, as in the previous case $H_4$ saturates $D_4$, as after adding an edge between $v_i$ and $v_j$ labeled with zero, the graph induced by these two vertices and $u_1$ and $u_2$ is isomorphic to $D_4$. We can add the edge $u_1u_2$ with the largest label, this does not affect the linearity of $sat_m(n, D_4)$. Thus, we are done.
	
	Finally, we treat the case $D_5$, for which the following labeling of $H$ gives $H_5$:
	\[ l(u_1v_i) = i \text{ and } l(u_2v_i) = 2n-3-i \text{, for all } i \in \{1,2, ..., n-2\}. \]
	This can be checked analogously to the previous two cases.
\end{proof}

The next statement is especially useful, as the way the host graph is defined in its proof will be reused multiple times later.

\begin{claim} \label{disjoint_neighborhood}
	Let $G$ be an edge-ordered graph, let $ab$ be the minimal edge in $G$, and assume that $deg(a) \ge 2$. If $N_G(a) \cap N_G(b) = \emptyset$, then $ssat_m(n,G) = \Theta(n)$.
\end{claim}

\begin{proof}
	The lower bound follows from Claim \ref{iso-edge-ssat}. For the upper bound, let $N_G(a) \setminus \{b\} =\{u_1, u_2, ..., u_s\}$ and $N_G(b) \setminus \{a\} = \{v_1, v_2, ..., v_t\}$, and assume that $G$ has $r$ vertices. Furthermore, let $H'$ be the graph obtained by taking the union of $G \setminus \{a, b\}$ and $\Bar{K}_{n-r+2}$ (the complement of the complete graph, that is, an independent set), whose set of vertices we denote by $\{w_1, w_2, ..., w_{n-r+2}\}$. Then, the host graph $H$ on $n$ vertices is obtained by adding edges $u_iw_k$ and $v_jw_k$ to $H'$, for $i \in [s]$, $j \in [t]$, and $k \in [n-r+2]$. We label them as follows: $L(u_iw_k) = l(u_ia).k$ and $L(v_jw_k) = l(v_jb).k$, while all the remaining edges keep the same labels as in $G$.
	
	Notice that $H$ semisaturates $G$ for the edges inside $\Bar{K}_{n-r+2}$. Indeed, if we add any edge of the form $w_iw_j$ for $i,j \in [n-r+2]$, it will have the role of the minimal edge $ab$ in the new copy of $G$. 
	
	We can greedily add a linear number of edges between $G \setminus \{a, b\}$ and the rest of the host graph $H$ to obtain a semisaturating host graph.
\end{proof}

A natural family that satisfies these conditions is the family of triangle-free graphs (such that $e_0$ is not isolated).

\begin{cor} \label{ssat_m triangle free}
	Let $G$ be a triangle-free edge-ordered graph, whose minimal edge is not isolated. Then $ssat_m(n, G) = \Theta(n)$.
\end{cor}

\begin{cor} \label{disjoint_nbhd_satm}
	Let $G$ be a non-bipartite edge-ordered graph. Let $e_0 = ab \in E(G)$ be the minimal edge and assume that $G-e_0$ is bipartite and $N_G(a) \cap N_G(b) = \emptyset$, then $sat_m(n,G) = \Theta(n)$.
\end{cor}

\begin{proof}
	We use the same construction as in Claim \ref{disjoint_neighborhood}. To verify that the host graph is $G$-free, notice that we obtain the host graph by taking the bipartite graph $G - ab$ and replacing the vertices $a$ and $b$ with an arbitrary number of vertices that are connected to $N_G(a) \cup N_G(b)$. This means that the host graph $H$ remains bipartite and consequently cannot contain a copy of the non-bipartite graph $G$.
\end{proof}

Corollary \ref{disjoint_nbhd_satm} implies that the cycles of odd length have a linear saturation function. We now discuss even cycles, as well as certain families of paths.

\begin{claim} \label{eocycles}
	Let $C_k$ be a cycle on $k \geq 5$ vertices with any edge-ordering. Then \[sat_m(n,C_k) = \Theta(n).\]
\end{claim}

\begin{figure}[ht]
	\centering
	\includegraphics[width=0.5\textwidth]{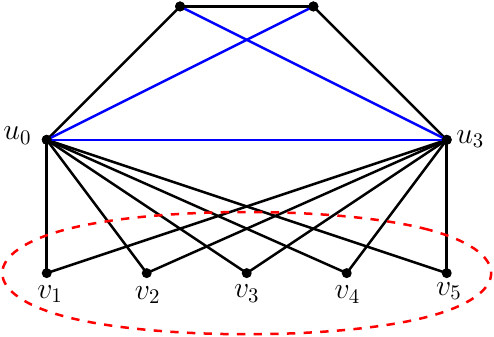}
	\caption{\label{fig: C6 host} A host graph which avoids the underlying graph $C_6$, and saturates the $C_6$ with any edge-ordering. The greedily added edges are denoted in blue, and the independent set $V$ of linear size in red.}
\end{figure}

\begin{proof}
	Let $H$ be the host graph obtained from the construction described in Claim \ref{disjoint_neighborhood}, and let $V = \{v_1, ..., v_{n-k+2}\}$ be the independent set in it. Let $u_1u_2$ be the minimal edge in $C_k$, and let $u_0$ (resp. $u_3$) be the other neighbor of $u_1$ (resp. $u_2$). Assume there is a copy of $C_k$ in $H$. Then clearly, there should be at least two vertices from $V$ in such a cycle, $v_i$ and $v_j$. Note that both of these vertices have degree two and  are both connected to $u_0$ and $u_3$. However, the longest cycle visiting the vertices $v_i, v_j, u_0, u_3$ has length $4$. Therefore, adding an edge $v_iv_j$ does not introduce a cycle of length $k \geq 5$. See Figure \ref{fig: C6 host} for the case $k = 6$.
\end{proof}

\begin{claim} \label{eopaths}
	Let $P_k$ be a path on the vertex set $\{u_1, u_2, ...,u_k\}$ for $k \geq 4$, such that $l(u_2u_3) \neq 0$ and $l(u_{k-2}u_{k-1}) \neq 0$, where $0$ is the minimal edge-label. Then $sat_m(n,P_k) = \Theta(n)$.
\end{claim}

\begin{proof}
	Following the same procedure as in Claim \ref{disjoint_neighborhood}, we only need to show that the host graph is $P_k$-free. We will show that its underlying graph does not have any path on $k$ vertices, regardless of its labels. 
	
	Consider first the case where $k \geq 6$ and $l(u_iu_{i+1}) = 0$, for some $i \in \{3,...,k-3\}$. Then, after we delete the vertices $u_i$ and $u_{i+1}$, we get two paths $u_1u_2...u_{i-1}$ and $u_{i+2}u_{i+3}...u_k$. Then we connect $u_{i-1}$ and $u_{i+2}$ with all the vertices in the set $V = \{v_1,..., v_{n-k+2}\}$ and obtain the host graph $H$. 
	
	Assume on the contrary that there is a copy of $P_k$ in $H$. Then there must be at least two vertices from $V$ in such a copy. Moreover, these two vertices must be connected either through $u_{i-1}$ or $u_{i+2}$. Suppose wlog. it is $u_{i+2}$, but then none of the vertices $u_j$, $j > i+2$ can be in the copy of the path $P_k$. However, this means that our path can have at most $i+2$ vertices, which is not enough since we assumed that $i \leq k-3$. Contradiction.
	
	It remains to check the case when $l(u_1u_2) = 0$ (or by symmetry $l(u_{k-1}u_k) = 0$). We can see that again, two vertices from the independent set $V$ are necessary, and they are connected through $u_3$. But this only gives a copy of $P_3$, while we assumed that $k \geq 4$. So, we again get a contradiction.
\end{proof}

Inspired by the construction provided in Claim \ref{D3-D5} that was used to prove the linearity of saturation function of $D_4$, we can prove the following semisaturation result.

\begin{claim} \label{linearfamilyD4}
	Let $G$ be an edge-ordered graph with labels $l: E(G) \rightarrow \mathbb{N}$ and let $ab$ be the minimal edge in $G$. Assume that $ab$ is not isolated and that for all $w \in N_G(a) \cap N_G(b)$ we have $l(aw) = l(bw)+1$. Then, $ssat_m(n,G) = \Theta(n)$.
\end{claim}

\begin{proof}
	Assume that $G$ has $k+2$ vertices for some $k \in \mathbb{N}$. Then we create $H$ starting with a copy of $G \setminus \{a,b\}$ keeping the same labels and $n-k$ isolated vertices denoted by $v_1, ...,v_{n-k}$. Then we add an edge between $v_i$ and $w$ for all $1 \leq i \leq n-k$ and $w \in N_G(a) \cup N_G(b)$ and label them with 
	\begin{equation*}
		L(v_iw) =
		\begin{cases}
			l(aw).i, & \text{if } w \in N_G(a)
			\\
			l(bw).i, & \text{if } w \in N_G(b) \setminus N_G(a).
		\end{cases}
	\end{equation*}
	
	As usual, it is simple to check that $H$ semisaturates $G$ for the edges inside $\Bar{K}_{n-k}$ and then we can greedily add at most a linear number of edges to get a semisaturating host graph.
\end{proof}

Claim \ref{linearfamilyD4} serves as a warm-up for the following, more general statement, for which we first need to introduce some notations.

\begin{defi}
	Let $l(E') = \{l(e) \colon e \in E'\}$, where $E' \subseteq E(G)$. We say that $l(E')$ forms an \emph{interval} in $l(G)$ if there is no other edge $e \in E(G)$ such that $l(e_1) < l(e) < l(e_2)$, for some $e_1, e_2 \in E'$.
\end{defi}

\begin{defi}
	Let $G$ be an edge-ordered graph with a labeling function $l: E(G) \rightarrow \mathbb{N}$ and $e_0 = ab$ its minimal edge. We denote by $\{u_1, u_2, ..., u_k\} := N_G(a) \cap N_G(b)$ and write $l(au_i) = a_i$ and $l(bu_i) = b_i$ for all $1 \leq i \leq k$. We assume that $a_1 < a_2 < ... < a_k$ and $b_1 < b_2 < ... < b_k$. Let \[A_1 = \{a_1, ..., a_{c_1}\}, A_2 = \{a_{c_1+1}, ..., a_{c_2}\},..., A_d = \{a_{c_{d-1}+1}, ..., a_{c_d}\},\]
	\[B_1 = \{b_1, ..., b_{c_1}\}, B_2 = \{b_{c_1+1}, ..., b_{c_2}\},..., B_d = \{b_{c_{d-1}+1}, ..., b_{c_d}\},\]
	for some $0 = c_0 \leq c_1 \leq ... \leq c_d = k$. 
	
	If for all $1 \leq i \leq d$, we have $A_i < B_i$ and $A_i \cup B_i$ forming an interval in $l(G)$, then $G$ is called a \emph{minimally interlacing} edge-ordered graph. 
\end{defi}

\begin{figure}[ht]
	\centering
	\includegraphics[width=0.5\textwidth]{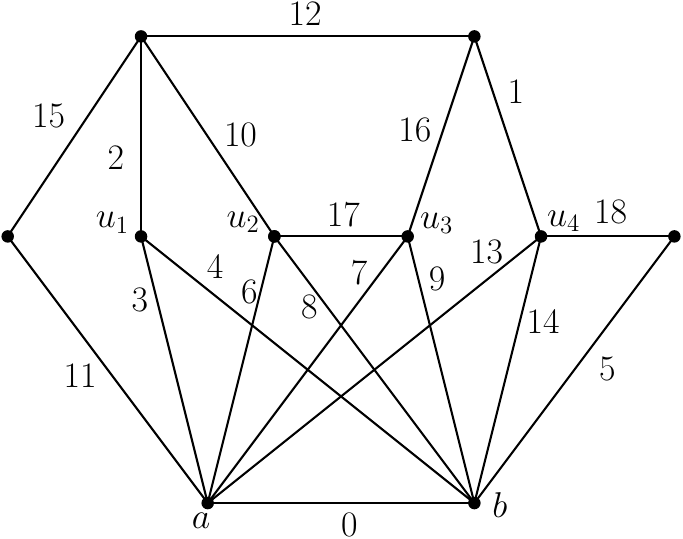}
	\caption{\label{fig: A1A2A3} $\{a_1, b_1\} = \{3,4\}, \{a_2, a_3, b_2, b_3\} = \{6,7,8,9\}, \{a_4, b_4\}= \{13, 14\}$.}
\end{figure}

See Figure \ref{fig: A1A2A3} for an example. 
The minimally interlacing edge-ordered graphs are labeled such that the vertices $a$ and $b$ in the minimal edge do not have a disjoint neighborhood, yet they also do not satisfy the assumptions of Theorem \ref{D_0 general} that would yield a superlinear semisaturation function. A natural question that arises is whether, any other edge-ordered graphs beyond those covered in Theorem \ref{D_0 general}, have a superlinear saturation function. In Remark \ref{increasing_order}, we discussed labelings of graphs for which Theorem \ref{D_0 general} does not apply. Even though minimally interlacing graphs represent only special cases of these families, the approach below might suggest that many of them indeed have a linear semisaturation function.

\begin{claim}
	Let $G$ be a minimally interlacing edge-ordered graph, then \[ssat_m(n,G) = \Theta(n).\]
\end{claim}
\begin{proof}
	Assume that $G$ has $m+2$ vertices for some $m \in \mathbb{N}$. We create $H$ from a copy of $G \setminus \{a,b\}$, keeping the same labels and from $n-m$ isolated vertices denoted by $V = \{v_1, ...,v_{n-m}\}$. Then we add an edge $v_iw$ for all $1 \leq i \leq n-m$ and $w \in N_G(a) \cup N_G(b)$ and label them as follows:
	\begin{equation*}
		L(v_iw) =
		\begin{cases}
			a_{c_s+1}.i.j, & \text{if } w = u_j, c_s+1 \leq j \leq c_{s+1},
			\\
			l(aw).i, & \text{if } w \in N_G(a) \setminus N_G(b),
			\\
			l(bw).i, & \text{if } w \in N_G(b) \setminus N_G(a).
		\end{cases}
	\end{equation*}
	Adding an edge in $V$ creates a copy of $G$ with some labeling. It remains to check that this copy is labeled as in $G$. We add an edge $v_iv_{i'}$ to $H$ for $1 \leq i < i' \leq n-m$. The edge $v_iv_{i'}$ has to play the role of $ab$ in the copy of $G$. It is enough to check that the neighbors of $v_i$ and $v_{i'}$ are labeled correctly. If $w \in N_G(a) \triangle N_G(b)$, then $l(aw).i$ is mapped to $l(aw)$ (resp. $l(bw).i$ is mapped to $l(bw)$). If $w \in N_G(a) \cap N_G(b)$, then for all $0 \leq s \leq d-1$ the interval $\{a_{c_s+1}.i.(c_s+1), ..., a_{c_s+1}.i.c_{s+1}, a_{c_s+1}.i'.(c_s+1),..., a_{c_s+1}.i'.c_{s+1}\}$ is mapped to the interval $\{a_{c_s+1}, ..., a_{c_{s+1}},b_{c_s+1}, ..., b_{c_{s+1}}\}$. For example, the interval $\{a_2, a_3, b_2, b_3\}$ in Figure \ref{fig: A1A2A3} is given in the host graph by the interval $\{6.i.2, 6.i.3, 6.i'.2, 6.i'.3\}$.
	
	We conclude that $H$ indeed semisaturates $G$ for a non-edge in the independent set $V$.
	Finally, we can add linearly many edges to $H$ to get a semisaturating host graph.
\end{proof}

\noindent \textbf{Extremal functions of edge-ordered graphs.} As noted in Remark \ref{rem:sme}, the extremal function $ex_e(n,G)$ provides an upper bound for all variants of the saturation function. Kucheriya and Tardos~\cite{kucheriyatardos2} characterized connected edge-ordered graphs with a linear extremal function. Note that any connected graph has at least a linear saturation function by Theorem \ref{isolated_edge}. This means that we can apply their results to get more examples of graphs with a linear saturation function. In particular, their results include trees with certain edge-orderings, which we did not discuss here. By Corollary \ref{ssat_m triangle free}, any edge-ordered tree (except a single edge) has a linear semisaturation function. However, we do not have the same result for saturation. 

On the other hand, the extremal functions for paths already have more complicated behaviour. It was shown in~\cite{gerbner} that the extremal functions of certain paths on four edges can be $\Theta(n\log n)$, or even $\binom{n}{2}$, while we showed that for almost all paths, the saturation function is linear. By Corollary \ref{ssat_m triangle free}, every edge-ordered forest has a linear semisaturation function, while the analogous result for saturation functions is unknown. It was shown in~\cite{gerbner} that for any edge-ordered star forest, there is an almost linear upper bound $ex(n,G) \leq n2^{(\alpha(n))^c}$, where $\alpha(n)$ denotes the inverse Ackermann function.

\section{Bounded saturation functions}

In this section, we only deal with $sat_m$ and $ssat_m$, the other variations are discussed in later sections. Claim \ref{iso-edge-ssat} characterizes edge-ordered graphs with bounded $ssat_m$, they are exactly those in which $e_0$ is isolated. In fact, it may be true that the same holds for $sat_m$. We consider special cases of edge-ordered graphs with isolated $e_0$ and prove that they have bounded $sat_m$.

As we showed in Theorem \ref{isolated_edge}, the saturation function can be bounded only if the minimal edge of a graph $G$ is isolated. In the following, we consider a construction which will help us prove that for many edge-ordered graphs the saturation function is indeed bounded.

\begin{defi}
	Let $G$ be an edge-ordered graph whose edges are labeled with the set $[m]$. Take three copies of $G$, which are denoted by $G_1$, $G_2$, and $G_3$, such that each $G_i$ is labeled with a map $k \rightarrow 3(k-1)+i$, for $k \in [m]$ and $i \in \{1,2,3\}$. We denote by $a_ib_i$ the corresponding minimal edge of $G_i$. Finally, we obtain $T(G)$ by merging the following pairs of vertices: $a_1$ and $b_3$, $a_2$ and $b_1$, $a_3$ and $b_2$.
\end{defi}

\begin{lem}
	Let $G$ be an edge-ordered graph labeled with the set $[m]$, and let $H$ be the union of $T(G)$ and a set of isolated vertices. Then $H$ semisaturates $e_0+G$ for any non-edge $uv$, where at least one of the vertices is isolated.
\end{lem}
\begin{proof}
	Assume that $u$ is an isolated vertex and $v$ is arbitrary. Then there must exist $i \in \{1,2,3\}$ such that $v \notin V(G_i)$. Therefore, $uv + G_i$ contains a copy of $e_0 + G$.
\end{proof}

Therefore, to show that $sat_m(n,e_0+G)$ is bounded for some $G$, it is enough to verify that $T(G)$ does not contain $e_0 + G$.

\begin{defi}
	We say that $G$ is an $A$-$B$ graph if it is obtained from two disjoint connected graphs $G_A$ and $G_B$ and a minimal edge $e_1 = ab$ such that $a \in V(G_A)$ and $b \in V(G_B)$.
\end{defi}

\begin{figure}[ht]
	\centering
	\includegraphics[width=0.5\textwidth]{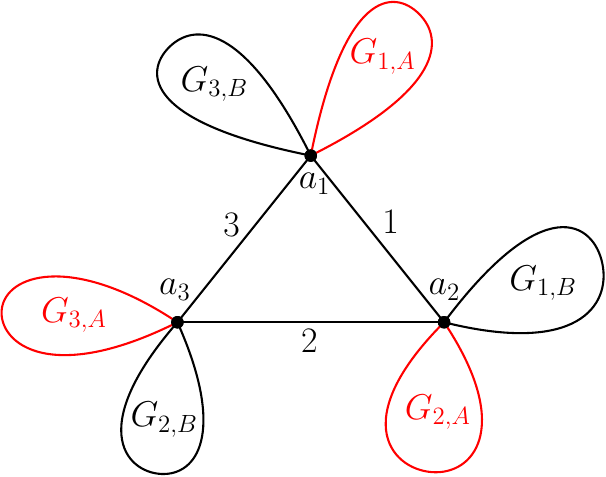}
	\caption{\label{fig:GA-GB} Graph $T(G)$ obtained by merging three copies of an $A$-$B$ graph $G$.}
\end{figure}

\begin{claim} \label{G_A - G_B}
	Let $F = e_0 + G$, where $G$ is an $A$-$B$ graph labeled with the set $[m]$. Then $sat_m(n,F) = O(1)$.
\end{claim}

\begin{proof}
	We consider a host graph $H$ which is a union of $T(G)$ and isolated vertices. We also assume that $a_ib_i$ connects the components $G_{i,A}$ and $G_{i,B}$, as shown in Figure \ref{fig:GA-GB}. 
	
	To show that $H$ is $F$-free, we assume the contrary. First, we consider the case when none of the three minimal edges of $H$ (that is, the edges $a_1a_2, a_2a_3, a_3a_1$) is contained in the copy of $F$. Notice that by deleting these edges from $H$ we would get a graph with three connected components such that each of them has $|V(G)| - 1$ vertices. Therefore, $H$ cannot contain $F$ in this case. Otherwise, one such edge, wlog. $a_1a_2$ is contained in the copy of $F$ and has to play the role of the minimal edge $e_0$ in the copy of $F$. However,  $H \setminus \{a_1, a_2\}$ has connected components only of size at most $|V(G)|-1$ and thus it cannot contain $G$, which gives a contradiction.
\end{proof}

We can easily see that any edge-ordered tree is an $A$-$B$ graph.

\begin{cor} \label{edgord_tree}
	Let $T$ be an edge-ordered tree, then $sat_m(n, e_0 + T) = O(1)$.
\end{cor}

\begin{obs}\label{G+F obs}
	Let $H(F)$ be a saturation graph of $F$ such that $sat_m(n,e_0 + F) = O(1)$. Let $G$ be an edge-ordered graph, such that\footnote{Here $l_{min}(G)$ denotes the minimal label of $G$.} $l_{min}(G) < l_{min}(F)$. Then $H(F)$ avoids $G+F$.
\end{obs}

Using this observation, we can extend Claim \ref{G_A - G_B} to certain disconnected graphs.

\begin{claim}
	Let $G = \sum_{i=1}^{k} G_i$ be a disjoint union of $A$-$B$ graphs such that, $|V(G_i)|\leq|V(G_{i+1})|$ and $l_{min}(G_i) < l_{min}(G_{i+1})$ for all $1 \leq i \leq k-1$. Then it follows that $sat_m(n,e_0 + G) = O(1)$.
\end{claim}

\begin{proof}
	Let $H$ be the union of a set of isolated vertices and of $\sum_{i=1}^{k} T(G_i)$ such that for edges $e_i\in T(G_i)$ and $e_j\in T(G_j)$ we have $l(e_i)<l(e_j)$. It is simple to check that $H$ indeed semisaturates $e_0+G$ apart from $O(1)$ edges. To show that it also avoids $e_0+G$, we assume the contrary and proceed by induction. The case $k=1$ follows by Claim \ref{G_A - G_B}. Then we assume that the claim holds for all $1 \leq j \leq k-1$. 
	
	We distinguish two cases. First, assume there is a copy $e_0+G$ in $H$ such that $G_k$ is in $T(G_k)$. In that case, by the inductive hypothesis either $e_0$ or  another $G_j$, for $j<k$ is also in $T(G_k)$, but this is not possible by Observation \ref{G+F obs}. Second, assume there is a copy of $G_k$ in $T(G_j)$ for some $j<k$. Since we assumed that $|V(G_j)|\leq|V(G_k)|$, it follows that the copy of $G_k$ contains at least one of the edges in $T(G_j)$ that corresponds to the minimal edge of $G_j$, since by deleting these three edges we get three connected components of size $|V(G_j)|-1$. This implies that $\sum_{i=1}^{j-1} T(G_i)$ contains a copy of $e_0 + \sum_{i=1}^{j-1} G_i$, which is not possible by the inductive hypothesis.
\end{proof}

We identify several other families of edge-ordered graphs for which $T(G)$ allows us to show that their saturation functions are bounded.

\begin{claim}
	Let $F = e_0 + G$, where $G$ is connected and satisfies the following property: for every pair of vertices $u,v \in V(G)$ (possibly $u=v$), the graphs $G \setminus \{u\}$, $G \setminus \{u,v\}$ are connected. Then $sat_m(n,F) = O(1)$.
\end{claim}
\begin{proof}
	We use again the same construction $T(G)$ as a host graph and we need to show that it avoids $F$. Notice first that a copy of $G$ cannot be identical to any of the graphs $G_1$, $G_2$ and $G_3$. So, it means that the copy of $G$ has to be found in the union of either two or three of these graphs. Also a copy of $G$ cannot contain any of the three edges $a_1a_2, a_2a_3, a_3a_1$, as there is no edge in that case which has the role of $e_0$.
	
	Assume wlog. that $G$ is in the union of $G_1$ and $G_2$, then $a_3$ is certainly in the copy of $G$, otherwise the graph would not be connected. But then $a_3$ is a vertex separator, which is not possible by assumption. We then assume that $G$ can be contained in all three copies. In that case, all three vertices $a_1, a_2, a_3$ have to be in the copy of $G$. Moreover, since $G$ is connected there should be at least two vertices among them, wlog. $a_1$ and $a_2$ such that there is a path of length at least two between them (recall that the edge $a_1a_2$ is not in the copy of $G$). So, if we delete the vertices $a_1$ and $a_2$, the copy of $G$ becomes disconnected which contradicts our assumption.
\end{proof}

\begin{cor} \label{e_0 complete}
	Let $K_r$ be a complete graph for $r \geq 2$ with arbitrary labeling. Then $sat_m(n,e_0 + K_r) = O(1)$.
\end{cor}

\begin{claim}
	Let $F = e_0 + G$, where $G$ is a connected graph whose minimal edge is $ab$ and $deg(a) \geq deg(b) > deg(v)$ for all $v \in V(G) \setminus \{a,b\}.$ Then $sat_m(n,F) = O(1)$.
\end{claim}
\begin{proof}
	We again want to show that $T(G)$ is $F$-free. Notice that the only vertices that have large enough degrees to have the role of the vertices $a$ and $b$, are the vertices $a_1$,$a_2$ and $a_3$. So, the minimal edge of $G$ must connect two of these vertices, but then there is no edge in $T(G)$ that can play the role of $e_0$, which leads to a contradiction.
\end{proof}

We now identify another family with a bounded saturation function, but this time without using the construction $T(G)$. Instead of three copies of the graph $G$, we now use only two and add a few additional edges.

\begin{thm} \label{G_peak}
	Let $F = e_0 + G$, such that there exists a vertex $v \in V(G)$, which we call the \emph{peak}, of degree $d \notin \{0,2\}$ whose incident edges are labeled with the set $\{1,...,d\}$ (the rest of the edges get larger labels) and $G-v$ is connected. Then \[sat_m(n,F) = O(1).\]
\end{thm}

\begin{proof}
	We proceed by constructing a host graph $H$ that semisaturates $F$. We first join two copies of $G$ at the peak $v$ and we denote the two disjoint parts by $G'$ and $G''$. Then we consider the neighbors $\{v_1, v_2, ...,v_d\} \subset V(G')$ of $v$ and we label each edge $l(vv_i) = 3(i-1)+1$ for $1 \leq i \leq d$. Additionally, we label by $l(vu_i) = 3(i-1)+2$, where $u_i$ is the neighbor of $v$ in $G''$ for all $1 \leq i \leq d$. Moreover, we add the edges between $u_1$ and all the $v_i$ and we label them with $l(u_1v_i) = 3i$. The remaining edges in $G'$ and $G''$ can be any labels greater than $3d$ that respect the edge-order inherited from $G$. See Figure \ref{fig: H for G 1-d} for a construction. As usual, the remaining vertices are isolated.
	
	\begin{figure}[ht]
		\centering
		\includegraphics[width=0.6\textwidth]{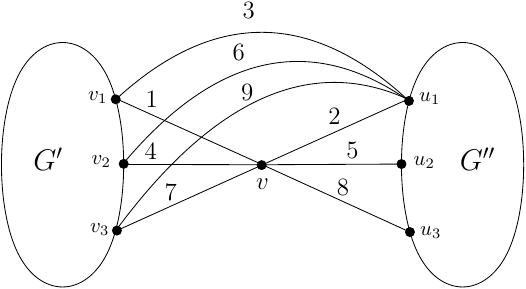}
		\caption{\label{fig: H for G 1-d} A saturation graph if the peak has degree $d=3$.}
	\end{figure}
	
	First, we want to show that $H$ is $F$-free. We treat the case $d=1$ separately by observing that no matter which edge we take to play the role of $e_0$, the copy of $G$ must be contained in either $G'$ or $G''$. That is not possible, since $G$ is larger than either of them. From now on, we consider the case $d \geq 3$. Assume on the contrary that $H$ is not $F$-free and assume that in a copy of $F$ we have a peak $p \notin N_H[v]$. In that case, we can notice that none of the edges connecting two vertices of $N_H[v]$ can be in the copy of $G$. However the rest of the graph $H$ contains two connected components with less than $|V(G)|$ vertices, a contradiction. 
	
	Let $p=v$ be the peak, in that case one of the edges $u_1v_1,...,u_1v_d$ must play the role of $e_0$ in the copy of $F$, because none of the edges incident to $v$ can be minimal, and the edges within $G'$ and $G''$ have all larger labels than the edges incident to $v$. Thus, if we remove $v$ from the copy of $G$, we get a disconnected graph, as some of its neighbors are in $G'$ and the other are in $G''$ contrary to our assumption that $G-v$ is connected. 	
	
	Assume now that the peak is $p=v_i$ for some $1 \leq i \leq d$. As $G$ is larger than $G'$, the copy of $G$ has a vertex outside $G'$. The peak $p$ must have at least one neighbor in $G'$ (as it has only two further incident edges in $H$) and such an edge has larger label than any edge  between $G'$, $p$, and $G''$. Therefore, no edge in the copy of $G-v$ can use the edges between $G'$, $p$, and $G''$ and so $G-v$ is a disconnected graph, a contradiction. A similar argument works for $p \in \{u_2, ..., u_d\}$. 	
	It remains to check if $u_1$ can be the peak. Clearly, $u_1v$ and $u_1v_1$ cannot be in the copy of $G$ as there would be no choice for the minimal edge. Thus, one neighbor of the peak has to be in $G''$, and for the same reason as in the previous case, we get that $u_1$ cannot be the peak.
	
	It remains to show that $H$ is $F$-saturated for the non-edges incident to isolated vertices. Indeed, if we take an isolated vertex and connect it with any vertex in $G'$ or $G''$, we get a new copy of $F$. If we connect it with $v$, then $u_1$ can have the role of the peak and together with $G'$ we get a copy of $F$.
	
	Finally, we can add a bounded number of edges greedily between non-isolated vertices to get a required host graph.		
\end{proof}

\begin{lem}
	Let $F = \sum_{i=1}^{k} F_i$ be such that $l(F_i) < l(F_j)$ for all $1 \leq i<j \leq k$ and $sat_m(n,e_0 + F_i) = O(1)$, for all $1 \leq i \leq k$, where $F_i$ are connected. Then it follows that $sat_m(n,e_0 + F) = O(1)$.
\end{lem}

\begin{proof}
	Let $H$ be a union of $\sum_{i=1}^{k} H(F_i)$ and isolated vertices, where $H(F_i)$ is the non-isolated part of the saturation graph of $F_i$. Additionally we label $H$ such that $l(H(F_i)) < l(H(F_j))$ for all $1 \leq i<j \leq k$. It is simple to see that $H$ semisaturates $e_0 + F$ apart from $O(1)$ edges. So, it remains to check that it also avoids $e_0 + F$. 
	
	We proceed by induction on $k$. The case $k = 1$ follows directly by assumption. Assume on the contrary that for $k$ the claim fails.
	As the claim holds for $k-1$ by inductive hypothesis, $\sum_{i=1}^{k-1} H(F_i)$ contains no copy of $e_0 + \sum_{i=1}^{k-1} F_i$. This implies that there is a copy of $F_j$ in $H(F_k)$ for some $j \leq k-1$, and clearly a copy of $F_k$ in $H(F_k)$ since $l(F_j) < l(F_k)$. But this is not possible by Observation \ref{G+F obs}.
\end{proof}

The previous result allows us to take a monotone union of any connected graphs that we previously showed to have bounded saturation function. Among others, this generalizes Corollary \ref{edgord_tree}.

\begin{defi}
	Let $F$ be a \emph{monotone forest} if $F = \sum_{i=1}^{k} T_i$ such that every $T_i$ is a tree and $l(T_i) < l(T_j)$ for all $1 \leq i<j \leq k$. 
\end{defi}

\begin{cor} \label{mon forests}
	For every monotone forest $F$, we have $sat_m(n,e_0 + F) = O(1)$.
\end{cor}

\section{\texorpdfstring{$sat_s$}{sat?} and \texorpdfstring{$ssat_s$}{ssat?} functions} \label{sats section}

In Definition \ref{definitionsAMS}, we stated two more definitions of saturation function. The reason why we do not have a unique definition is that we have a choice on how to label the missing edge. So far, we worked with the definition where a new edge always gets the minimal label. Recall that, $sat_s(n,G)$ denotes the minimum number of edges in an edge-ordered graph $H$ on $n$ vertices in which for any new edge there exists some label (inserted somewhere in the linear order of edges) which introduces a copy of $G$.

As mentioned in Remark \ref{rem:sme}, we have $sat_s(n,G) \leq sat_m(n,G)$, which implies that all the upper bounds that we proved for $sat_m$ definition hold for $sat_s$ as well. This also includes the general $O(n \log n)$ upper bound for semisaturation. Using an analogous argument as in the proof of Theorem \ref{D_0 general}, where the assumption for the minimal edge is replaced by the assumptions on every edge of the graph, we conclude that no dichotomy arises in this framework either. Therefore, the result for this definition holds for a more restricted family of graphs. As the proof follows the same lines as the proof of Theorem \ref{D_0 general}, we omit it here. 

\begin{claim}\label{claim:sats}
	Let $G$ be an edge-ordered graph. If for every edge $ab \in E(G)$ there exist $v, w \in N_G(a) \cap N_G(b)$ such that $l(av) < l(aw)$ and $ l(bv) > l(bw)$, then $sat_s(n,G) \geq ssat_s(n,G)=\Omega(n\sqrt{\log n})$. 
\end{claim}

Observe that none of the diamond graphs shown in Figure \ref{fig:D012345} satisfy the assumptions of Claim \ref{claim:sats}. However, we show that a complete graph on at least $4$ vertices can always be labeled such that it satisfies these assumptions.

\begin{claim}\label{claim:satsKk}
	For every $d \ge 4$, there exists an edge-ordered graph $G_k$ with underlying graph $K_k$, such that $sat_s(n,G_k) \geq ssat_s(n, G_k)=\Omega(n\sqrt{\log n})$.
\end{claim}

\begin{figure}[ht]
	\centering
	\includegraphics[width=0.35\textwidth]{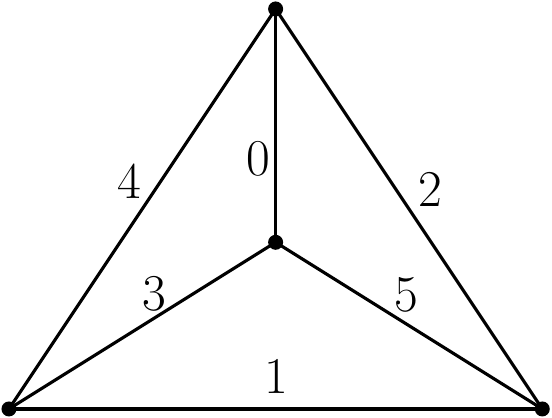}
	\caption{\label{fig:K4 sats} An example of edge-ordered $K_4$ with a superlinear $sat_s$ function.}
\end{figure}

\begin{proof}
	If $d=4$, we can take a copy of $K_4$ and label its edges as shown in Figure \ref{fig:K4 sats}. Moreover, for $K_d$ with $d \ge 5$, we proceed as follows: label the vertices with $[d]$ and for each pair $i \ne j$ in $[d]$, we assign to the edge $ij$ the label \[\min(|j-i|, d-|j-i|),\] which represents the cyclic distance between $i$ and $j$ modulo $d$. Then we perturb the labels slightly to ensure that all edge labels are distinct. It is easy to check that this labeling satisfies the assumptions of Corollary \ref{claim:sats} for every $d\ge 5$.
\end{proof}

Note that, by the inequality $sat_s(n,G) \leq sat_m(n,G)$ and Theorem \ref{weak dichotomy sats}, any graph $G$ with a linear $sat_m$ function and contains no isolated edges also has a linear $sat_s$ function. Moreover, the approach used in Claim \ref{disjoint_neighborhood} can be extended to identify an even larger family of edge-ordered graphs with linear semisaturation functions. 

\begin{claim} \label{disj_nbhd_ssats}
	Let $G$ be an edge-ordered graph with no isolated edges. If there exists an edge $uv \in E(G)$ such that $N_G(u) \cap N_G(v) = \emptyset$, then $ssat_s(n,G) = \Theta(n)$.
\end{claim}

For the proof, we refer the reader to Claim \ref{disjoint_neighborhood} for the construction of a host graph with a linear number of edges. The same construction applies here, with the condition on the minimal edge replaced by an analogous condition on an arbitrary edge of $G$.

\section{\texorpdfstring{$sat_e$}{sat?} and \texorpdfstring{$ssat_e$}{ssat?} functions}

The most natural saturation definition for edge-ordered graphs is $sat_e(n,G)$, where we allow the new edge to have any label inserted in the linear order of the edges of the host graph $H$. In this section, we state several results about this variant. We left these to be the last as most of the forthcoming results follow from the results that we have already proved.

Recall that for all edge-ordered graphs we have $sat_m(n,G) \leq sat_e(n,G)$. Therefore, all the lower bounds that we showed for $sat_m$ also hold for $sat_e$ definition. Most importantly, this gives us a family of edge-ordered graphs with superlinear semisaturation (and saturation) function. The next result follows directly from Theorem \ref{D_0 general}.

\begin{cor} \label{D_0 general sat_e}
	Let $G$ be an edge-ordered graph. If $ab \in \{e_0, e_{\max}\}$ and there exist two vertices $v,w \in N_G(a) \cap N_G(b)$ such that $l(av) < l(aw)$ and $ l(bv) > l(bw)$, then \[sat_e(n,G) \geq ssat_e(n,G)=\Omega(n\sqrt{\log n}).\]
\end{cor}

We now show the upper bound for the diamond graph $D_0$, which can be seen in Figure \ref{fig:D012345}. We make use of Theorem \ref{D_0 upper bound} which was proven for the $sat_m$ definition.

\begin{thm} \label{D_0 sate upper bound}
	\[ssat_e(n,D_0) \leq sat_e(n,D_0)=O\left(n \frac{\log n}{\log \log n}\right).\]
\end{thm}

\begin{proof}
	We start with $H = K_{k,n}$, the host graph that we constructed in the proof of Theorem \ref{D_0 upper bound}. Recall that $H$ was labeled by $\mathbb{R_{>}}$. We now add two more vertices $u_1$ and $u_2$, and connect them with the independent set in $H$ that we denote by $V = \{v_1, \dots, v_n\}$. We label the edges as follows:
	\[l(u_1v_i) = -3.i \text{ and } l(u_2v_i) = -1.i \text{ for all } 1 \le i \le n.\] We also add the edge $u_1u_2$ and label it by $l(u_1u_2) = -2$. Notice that by adding any edge between two vertices in $V$ and labeling it with any positive number, we obtain a copy of $D_0$, where the new edge plays the role of $bv$, the maximal edge in $D_0$ (see Figure \ref{fig:D012345}). Recall that from Theorem \ref{D_0 upper bound}, if any edge in $V$ is added with a negative label or zero, a new copy of $D_0$ is created. Therefore, we get a host graph that semisaturates $D_0$ for any non-edge in $V$ with any real label.
	
	It remains to check that the host graph avoids $D_0$, so we assume the contrary. Notice that by deleting the edge $u_1u_2$, we get a complete bipartite graph. Therefore, since $D_0$ is not bipartite, $u_1u_2$ has to be in the copy of $D_0$. Even more, $u_1u_2$ has to play the role of the minimal edge $ab$, because when any non-minimal edge of $D_0$ is deleted, the obtained subgraph is non-bipartite. But $u_1u_2$ cannot be minimal, because $l(u_1v_i) < l(u_1u_2)$ for all $1 \le i \le n$. This implies that $u_1u_2$ cannot be the minimal edge, which is a contradiction. Finally, we can add edges greedily between the $k+2$ vertices without introducing a copy of $D_0$, which does not affect the order of magnitude of the saturation function.
\end{proof}

Recall that by Theorem \ref{e_0e_max-isolated} if $G$ is an edge-ordered graph in which $e_0$ and $e_{\max}$ are isolated, then either $sat_e(n,G) = O(1)$ or $sat_e(n,G) = \Omega(n)$. Otherwise, $sat_e(n,G) = \Omega(n)$. Here we show why $\Omega(n)$ was necessary in the statement already for the $e_0+G+e_{\max}$ part, a behaviour that is different from that of $sat_m$ function in Theorem \ref{isolated_edge}.

\begin{thm} \label{complete_superlinear}
	For every $k\ge 5$, there exists an edge-ordered graph $G_k$ with underlying graph $K_k$, such that \[sat_e(n,e_0+G_k+e_{\max})=\Omega (n\sqrt{\log n}).\]
\end{thm}

\begin{proof}
	Let $G_k$ be an edge-ordered graph with underlying graph $K_k$. We denote $F_k=e_0+G_k+e_{\max}$.
	Let $H$ be a host graph saturating $F_k$ for $sat_e$. Take a non-edge. Adding it as a new edge $e$ with any label should introduce a copy of $F_k$. Note that there are only finitely many essentially different labels (depending on where the label fits in the linear order of labels in $F_k$), so we concentrate only on them. 
	
	We claim that there exists a label such that, in the introduced copy of $F_k$ the new edge $e$ does not play the role of $e_0$ or $e_{\max}$. Assume the contrary and take the smallest label such that $e$ plays the role of $e_{\max}$ in a copy $F'_k$ of $F_k$ (this exists, as when we choose the largest label for $e$ then $e$ must play the role of $e_{\max}$). Take the next smallest possible label for $e$ (this again must exist as the current label could not have been the smallest). Then $e$ must play the role of $e_0$ in another copy $F''_k$ of $F_k$. 
	
	Let $G'_k$ and $G''_k$ be the $K_k$ part of $F'_k$ and $F''_k$, respectively. Notice that all the labels of the edges of $G'_k$ are smaller than the labels of the edges of $G''_k$, except possibly for an edge appearing in both of them. Therefore, their vertex sets intersect in at most two vertices. Moreover, the edge that plays the role of $e_{\max}$ in $F''_k$ can share at most one vertex with $G'_k$. Therefore, there are at least two vertices in $G'_k$, which are not in $F''_k$. This edge can also play the role of $e_0$ in $F''_k$ instead of $e$, contradicting the assumption that $H$ avoids $F_k$.
	
	Thus, for each non-edge we can add it with an appropriate label such that we get a copy of $F_k$ in which the new edge $e$ plays the role of an edge in $G_k$. That is, adding $e$ with this label must introduce a new copy of $G_k$. This means that $H$ is semisaturating $G_k$ for $ssat_s$, therefore $sat_e(n,e_0+G_k+e_{\max})=|H|\ge ssat_s(n,G_k)$. Applying Claim \ref{claim:satsKk} gives that there is an appropriate $G_k$ such that $ssat_s(n,G_k) = \Omega(n\sqrt{\log n})$, finishing the proof.		
\end{proof}

\begin{figure}[ht]
	\centering
	\includegraphics[width=0.9\textwidth]{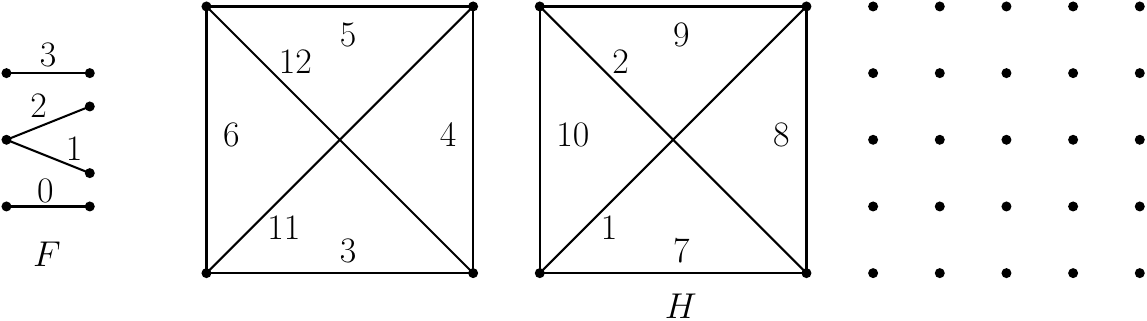}
	\caption{\label{fig:cherryhost}  $F=e_0+G+e_{\max}$, where $G$ is the cherry graph and an almost host graph $H$ for $F$.}
\end{figure}

One might wonder if the main observation in the previous proof holds more generally, that is, if $sat_e(n,e_0+G+e_{\max})=\Omega(ssat_s(n,G))$ always. We have seen that this holds if $G$ is an edge-ordered complete graph on at least $5$ vertices. Here we give an example showing that this is not true for a graph $G$ with two edges. 

Let $G$ be the unique edge-ordered cherry graph (two edges on three vertices). By Claim \ref{bounded-ssat_s}, we have $ssat_s(n,G)=\Omega (n)$. On the other hand, we claim that $F=e_0+G+e_{\max}$ has $sat_e(n,F)=O(1)$. Let $H$ be the graph on $n$ vertices such that $n-8$ of the vertices are isolated, and on the remaining $8$ vertices we have two $K_4$'s, labeled as shown in Figure \ref{fig:cherryhost}. It is easy to check that $H$ avoids $F$ yet adding any new edge incident to at least one isolated vertex with any label introduces a copy of $F$. Therefore, by greedily adding further edges between the two $K_4$'s if necessary, we obtain a saturating host graph on $n$ vertices and $O(1)$ edges.

Furthermore, notice that if two copies $G'$ and $G''$ of a graph $G$, such that $l(G') < l(G'')$, can be edge-disjointly added on $|V(G)|+1$ vertices, then $sat_e(n,e_0+G+e_{\max}) = O(n)$. Indeed, the host graph can be obtained by adding a vertex-disjoint cherry with minimal and maximal edge to this graph on $|V(G)|+1$ vertices (the two copies of $G$ have labels that form two intervals of $l(G)$), as well as arbitrary many isolated vertices. This graph clearly avoids $e_0+G+e_{\max}$ and by adding an edge with any label between two isolated vertices we get a copy of $e_0+G+e_{\max}$. Thus, by adding $O(n)$ edges, we get a host graph. With this construction, we can produce many examples that have at most linear saturation function.

\begin{figure}[ht]
	\centering
	\includegraphics[width=0.35\textwidth]{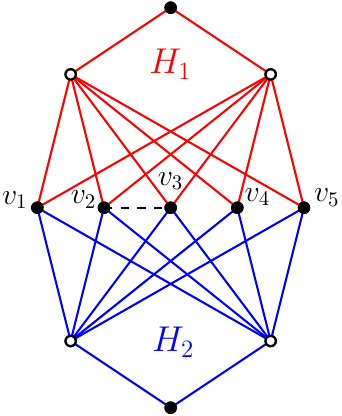}
	\caption{\label{fig:satC5}  The bipartite host graph that saturates $C_5$ for a non-edge in $\{v_1,v_2,v_3,v_4,v_5\}$, obtained by gluing two host graphs $H_1$ and $H_2$, such that $l(H_1) < l(H_2)$.}
\end{figure}

We can observe that the analogue of Corollary \ref{disjoint_nbhd_satm} also holds for the $sat_e$ definition. In comparison with Corollary \ref{disjoint_nbhd_satm}, the difference is that now we need to glue together two host graphs (one saturating for the minimal label and one for the maximal label), i.e., identify their subsets of vertices that form the `large' independent set, such that the resulting graph remains bipartite. For example, see Figure \ref{fig:satC5} for the resulting host graph when we saturate a cycle $C_5$. Therefore, we need to add a condition to ensure that $G - e_{\max}$ is also bipartite. Note that the host $H_1$ and $H_2$ graphs should be labeled such that $l(H_1) < l(H_2)$ in order to ensure that a new edge inside the large independent set can receive any label.

\begin{cor} \label{disjoint_nbhd_sate}
	Let $G$ be a non-bipartite edge-ordered graph. Let $e_0 = ab$ and $e_{\max} = cd$ be such that $G-e_0$ and $G-e_{\max}$ are bipartite and $N_G(a) \cap N_G(b) = N_G(c) \cap N_G(d) = \emptyset$. Then $sat_e(n,G) = \Theta(n)$.
\end{cor}

From this, we get the following explicit infinite family of graphs that have linear $ssat_e$ and $sat_e$ functions. 

\begin{cor} \label{eeocycles}
	Let $C_{2k+1}$ be an edge-ordered cycle. Then for every $k \geq 1$, \[sat_e(n,C_{2k+1}) = \Theta(n) \text{ and } ssat_e(n,C_{2k+1}) = \Theta(n).\]
\end{cor}

To finish this section, we also present an infinite family of edge-ordered graphs with a bounded $sat_e$ function. Since all edge-ordered matching graphs with a fixed number of edges are isomorphic to each other, in this case, the saturation problem is equivalent to the saturation problem of unordered graphs. Therefore, we get from~\cite{kaszonyituza} that:

\begin{cor} \label{matchings}
	Let $M_{k}$ be an edge-ordered matching. Then \[sat_e(n, M_{k}) = O(1), \text{ for every k } \geq 1.\]
\end{cor}

\section{Open questions}

We conclude by highlighting the most interesting problems left open.

We do not have any example of an edge-ordered graph with a linear $sat_m$ saturation function whose minimal edge is isolated. The construction $T(G)$ (see Figure \ref{fig:GA-GB}) appears to be a promising candidate to show that no such graph exists.

\begin{conj}
	Let $F = e_0 + G$, then $sat_m(n,F) = O(1)$.
\end{conj}

Note that the analogous statement for the $sat_e$ definition, where we replace $e_0$ with $e_0+e_{\max}$, is not true, as we have seen earlier.

Claim \ref{disjoint_neighborhood} shows that if the end vertices of the non-isolated minimal edge of $G$ have disjoint neighborhoods, then the semisaturation function of $G$ is linear. So far, we do not have a counterexample, and it is possible that this family also has a linear saturation function.

\begin{conj}
	Let $G$ be an edge-ordered graph and let $ab$ be the minimal edge in $G$ such that $deg(a) \ge 2$. If $N_G(a) \cap N_G(b) = \emptyset$, then $sat_m(n,G) = \Theta(n)$.
\end{conj}

The previous two conjectures raise the question whether in general for a given graph, the orders of magnitude of $sat_m$ and $ssat_m$ are always the same. We do not have a counterexample for this claim. However, they do not have necessarily the same exact values, for example one can show that $sat_m(n, M_k) = 3(k-1)$ and $ssat_m(n, M_k) = k+1$, for $k \ge 2$. 
For $sat_e$ definition, we already have some counterexamples. Theorem \ref{bounded-ssat_e} and Theorem \ref{complete_superlinear} give edge-ordered graphs with a bounded $ssat_e$ function, but a superlinear $sat_e$ function.

The main open question that remains is to find a general upper bound for the saturation function. We expect to have a non-trivial upper bound, and we state it in terms of $sat_e$, which implies the same  upper bound for $sat_m$.

\begin{conj}
	For every edge-ordered graph $G$, we have $sat_e(n,G) = O(n^{1+o(1)})$.
\end{conj}

It seems plausible to expect that the same upper bound as for semisaturation can be proved. Therefore, we also state the following stronger conjecture.

\begin{conj}
	For every edge-ordered graph $G$, we have $sat_e(n,G) = O(n \log n)$.
\end{conj}

Besides that, it might be possible that the general upper bound depends on the chosen definition. 

We can check that the construction given in the proof of Theorem \ref{nlogn upper bound}, indeed gives a saturation upper bound for many graphs. However, in some cases, such as certain families of trees, the construction only gives a semisaturation upper bound. 

The semisaturation problem can be studied separately and we managed to understand the behaviour of semisaturation functions slightly better. The major question is whether the upper bound can be improved, or at least can we find an example when the bound is achieved.

\begin{prob}
	Is there an edge-ordered graph $G$ such that $ssat_e(n,G) = \Theta(n \log n)$?
\end{prob}

Studying thoroughly the diamond graph $D_0$ was helpful to understand how the saturation functions behave for more general families of edge-ordered graphs. In fact, we were able to extend the lower and upper bounds to infinite families of graphs. One possible direction is to study further the saturation function of $D_0$, which would again most likely lead to a more general result. On the other hand, $D_0$ is the only connected edge-ordered graph for which we showed an upper bound $O\left(n \frac{\log n}{\log \log n}\right)$. Can we extend this result to a more general family of connected edge-ordered graphs?

We were able to show at least three different behaviours of saturation functions. A natural question is if there are more than that, or more generally, are there finitely or infinitely many different orders of magnitude for saturation functions of edge-ordered graphs? Analogous question for semisaturation functions is also open.

Recently, the saturation functions of sequences were studied~\cite{seqsat}. They showed that the saturation function is either bounded or at least linear. It would be interesting to investigate whether an example with a superlinear saturation function can be found, or to show that the usual dichotomy holds. Regarding semisaturation functions, they showed that they are at most linear, which is not the case for edge-ordered graphs.

\subsection*{Acknowledgements}

We would like to thank the anonymous referees for their comments and suggestions, which greatly helped us improve the presentation of the paper.

Part of this research was performed while the first author was visiting the Institute for Pure and Applied Mathematics (IPAM), which is supported by the National Science Foundation (Grant No. DMS-1925919).

Research supported by the J\'anos Bolyai Research Scholarship of the Hungarian Acade\-my of Sciences; by the National Research, Development and Innovation Office -- NKFIH under grants K 132696 and FK 132060; by the \'UNKP-23-5 New National Excellence Program of the Ministry for Innovation and Technology from the source of the National Research, Development and Innovation Fund; and by the ERC Advanced Grant ``ERMiD''. 
This research has been implemented with the support provided by the Ministry of Innovation and Technology of Hungary from the National Research, Development and Innovation Fund, financed under the  ELTE TKP 2021-NKTA-62 funding scheme, and by the EXCELLENCE-24 project no.~151504, Combinatorics and Geometry, of the NRDI Fund.
	
\bibliographystyle{plainurl}
\bibliography{edgesaturation}

@article {dong007,
	AUTHOR = {Dong, Jinquan and Liu, Yanpei},
	TITLE = {On the decomposition of graphs into complete bipartite graphs},
	JOURNAL = {Graphs Combin.},
	FJOURNAL = {Graphs and Combinatorics},
	VOLUME = {23},
	YEAR = {2007},
	NUMBER = {3},
	PAGES = {255--262},
	ISSN = {0911-0119},
	MRCLASS = {05C70},
	MRNUMBER = {2320580},
	MRREVIEWER = {R. Balakrishnan},
}

@article {boskesz023,
	AUTHOR = {Bo\v{s}kovi\'{c}, Vladimir and Keszegh, Bal\'{a}zs},
	TITLE = {Saturation of {O}rdered {G}raphs},
	JOURNAL = {SIAM J. Discrete Math.},
	FJOURNAL = {SIAM Journal on Discrete Mathematics},
	VOLUME = {37},
	YEAR = {2023},
	NUMBER = {2},
	PAGES = {1118--1141},
	ISSN = {0895-4801},
	MRCLASS = {05C35},
	MRNUMBER = {4602504},
}

@article{kaszonyituza,
	author = {K\'aszonyi, L\'aszl\'o and Tuza, Zsolt},
	title = {Saturated graphs with minimal number of edges},
	journal = {Journal of Graph Theory},
	volume = {10},
	number = {2},
	pages = {203-210},
	year = {1986}
}

@misc{alon,
      title={On bipartite coverings of graphs and multigraphs}, 
      author={Noga Alon},
      year={2023},
      eprint={2307.16784},
      archivePrefix={arXiv},
      primaryClass={math.CO}
}

@article {katona,
	AUTHOR = {Katona, Gyula and Szemer\'{e}di, Endre},
	TITLE = {On a problem of graph theory},
	JOURNAL = {Studia Sci. Math. Hungar.},
	FJOURNAL = {Studia Scientiarum Mathematicarum Hungarica. Combinatorics,
	Geometry and Topology (CoGeTo)},
	VOLUME = {2},
	YEAR = {1967},
	PAGES = {23--28},
	ISSN = {0081-6906,1588-2896},
	MRCLASS = {05.60},
	MRNUMBER = {210619},
	MRREVIEWER = {F.\ Harary},
}

@article{erdoshajnalmoon,
 author = {Paul Erd\H{o}s and Andr{\'a}s Hajnal and John W. Moon},
 journal = {The American Mathematical Monthly},
 number = {10},
 pages = {1107--1110},
 publisher = {Mathematical Association of America},
 title = {A Problem in Graph Theory},
 volume = {71},
 year = {1964}
}

@article{brualdicao,
title = {Pattern-avoiding (0,1)-matrices and bases of permutation matrices},
journal = {Discrete Applied Mathematics},
volume = {304},
pages = {196-211},
year = {2021},
issn = {0166-218X},
doiTemp = {https://doi.org/10.1016/j.dam.2021.07.039},
urlTemp = {https://www.sciencedirect.com/science/article/pii/S0166218X21003061},
author = {Richard A. Brualdi and Lei Cao}
}

@article{fulkesz,
  author    = {Radoslav Fulek and
               Bal{\'{a}}zs Keszegh},
  title     = {Saturation Problems about Forbidden 0-1 Submatrices},
  journal   = {{SIAM} J. Discret. Math.},
  volume    = {35},
  number    = {3},
  pages     = {1964--1977},
  year      = {2021},
  bibsource = {dblp computer science bibliography, https://dblp.org}
}

@article{geneson,
	author    = {Jesse Geneson},
	title     = {Almost all Permutation Matrices have Bounded Saturation Functions},
	journal   = {Electron. J. Comb.},
	volume    = {28},
	number    = {2},
	pages     = {P2.16},
	year      = {2021},
	urlTemp       = {https://www.combinatorics.org/ojs/index.php/eljc/article/view/v28i2p16},
	bibsource = {dblp computer science bibliography, https://dblp.org}
}

@article {berendsohn,
	AUTHOR = {Berendsohn, Benjamin Aram},
	TITLE = {An exact characterization of saturation for permutation
	matrices},
	JOURNAL = {Comb. Theory},
	FJOURNAL = {Combinatorial Theory},
	VOLUME = {3},
	YEAR = {2023},
	NUMBER = {1},
	PAGES = {Paper No. 17, 35},
	ISSN = {2766-1334},
	MRCLASS = {05D99 (05A05 05B20)},
	MRNUMBER = {4565304},
	MRREVIEWER = {Luis\ Boza},
}

@article {gerbner,
	AUTHOR = {Gerbner, D\'{a}niel and Methuku, Abhishek and Nagy, D\'{a}niel
	T. and P\'{a}lv\"{o}lgyi, D\"{o}m\"{o}t\"{o}r and Tardos,
	G\'{a}bor and Vizer, M\'{a}t\'{e}},
	TITLE = {Tur\'{a}n problems for edge-ordered graphs},
	JOURNAL = {J. Combin. Theory Ser. B},
	FJOURNAL = {Journal of Combinatorial Theory. Series B},
	VOLUME = {160},
	YEAR = {2023},
	PAGES = {66--113},
	ISSN = {0095-8956,1096-0902},
	MRCLASS = {05C35 (52C10)},
	MRNUMBER = {4530857},
	MRREVIEWER = {Jen\H{o}\ Lehel},
}

@article {kucheriyatardos1,
	AUTHOR = {Kucheriya, Gaurav and Tardos, G\'{a}bor},
	TITLE = {A characterization of edge-ordered graphs with almost linear
	extremal functions},
	JOURNAL = {Combinatorica},
	FJOURNAL = {Combinatorica. An International Journal on Combinatorics and
	the Theory of Computing},
	VOLUME = {43},
	YEAR = {2023},
	NUMBER = {6},
	PAGES = {1111--1123},
	ISSN = {0209-9683,1439-6912},
	MRCLASS = {05C35},
	MRNUMBER = {4670572},
}

@misc{kucheriyatardos2,
	title={On edge-ordered graphs with linear extremal functions}, 
	author={Kucheriya, Gaurav and Tardos, G\'{a}bor},
	year={2023},
	eprint={2309.10558},
	archivePrefix={arXiv},
	primaryClass={math.CO}
}

@article{graphsatsurvey,
	author = {Bryan L. Currie and
	Jill R. Faudree and
	Ralph J. Faudree and
	John R. Schmitt},
	title = {A Survey of Minimum Saturated Graphs},
	journal = {The Electronic Journal of Combinatorics},
	volume = {DS19},
	year = {2021},
}

@Article{curves,
author={Keszegh, Bal{\'a}zs
and P{\'a}lv{\"o}lgyi, D{\"o}m{\"o}t{\"o}r},
title={The Number of Tangencies Between Two Families of Curves},
journal={Combinatorica},
year={2023},
month={Oct},
day={01},
volume={43},
number={5},
pages={939-952}
}

@article{xmon,
title = {On the number of tangencies among 1-intersecting x-monotone curves},
journal = {European Journal of Combinatorics},
volume = {118},
pages = {103929},
year = {2024},
issn = {0195-6698},
doiTemp = {https://doi.org/10.1016/j.ejc.2024.103929},
URLTemp = {https://www.sciencedirect.com/science/article/pii/S0195669824000143},
author = {Eyal Ackerman and Bal\'azs Keszegh}
}

@article {seqsat,
	AUTHOR = {Anand and Geneson, Jesse and Kaustav, Suchir and Tsai,
	Shen-Fu},
	TITLE = {Sequence saturation},
	JOURNAL = {Discrete Appl. Math.},
	FJOURNAL = {Discrete Applied Mathematics. The Journal of Combinatorial
	Algorithms, Informatics and Computational Sciences},
	VOLUME = {360},
	YEAR = {2025},
	PAGES = {382--393},
	ISSN = {0166-218X,1872-6771},
	MRCLASS = {05A05},
	MRNUMBER = {4806441},
}

@article{Furedi1998MinimalOG,
	title={Minimal Oriented Graphs of Diameter 2},
	author={Zolt{\'a}n F{\"u}redi and Peter Hor{\'a}k and Chandra M. Pareek and Xuding Zhu},
	journal={Graphs and Combinatorics},
	year={1998},
	volume={14},
	pages={345-350},
}

@article{Pikhurko, 
	title={The Minimum Size of Saturated Hypergraphs}, 
	volume={8}, 
	DOITemp={10.1017/S0963548399003971}, 
	number={5}, 
	journal={Combinatorics, Probability and Computing}, 
	publisher={Cambridge University Press}, 
	author={Pikhurko, Oleg}, 
	year={1999}, 
	pages={483--492}
}

@misc{perscommalonfox,
  author = "Noga Alon and Jacob Fox",
  howpublished = "personal communication",
  year = "2024"  
}
	
\end{document}